# Effects of water currents on fish migration through a Feynman-type path integral approach under $\sqrt{8/3}$ Liouville-like quantum gravity surfaces

Paramahansa Pramanik[1] 

**Abstract**

A stochastic differential game theoretic model has been proposed to determine optimal behavior of a fish while migrating against water currents both in rivers and oceans. Then, a dynamic objective function is maximized subject to two stochastic dynamics, one represents its location and another its relative velocity against water currents. In relative velocity stochastic dynamics, a Cucker–Smale type stochastic differential equation is introduced under white noise. As the information regarding hydrodynamic environment is incomplete and imperfect, a Feynman type path integral under $\sqrt{8/3}$ Liouville-like quantum gravity surface has been introduced to obtain a Wick-rotated Schrödinger type equation to determine an optimal strategy of a fish during its migration. The advantage of having Feynman type path integral is that, it can be used in more generalized nonlinear stochastic differential equations where constructing a Hamiltonian–Jacobi–Bellman (HJB) equation is impossible. The mathematical analytic results show exact expression of an optimal strategy of a fish under imperfect information and uncertainty.

**JEL classification** Primary C73 · Secondary C61

## Introduction

Migration of fishes is an important factor for environment surrounding of it and also animals based on them. Most of the times fishes migrate for spawning and foraging for food. Young fishes usually leave their spawning areas and go to places where they become adults. On the other hand, adult fishes move to the spawning area and then return to the feeding ground. During migrations fishes travel a long distance, and furthermore, migrations of some adult fishes such as *Plecoglossus altivelis* (Ayu), *Oncorhynchus masou* (Yamame) and *Poecilia reticulata* (guppy) toward breeding grounds are against the water current (Yoshioka 2017; Yoshioka and Yaegashi 2018; Yoshioka et al. 2019) and fish like larger brown trout initiates the downstream migration for spawning (Jonsson and Jonsson 2002). As availability of food at a particular destination is probabilistic, a fish has to use all directions irrespective of water currents. Therefore, adult movements are directional rather than passive. This paper considers only migratory fish movements and effectiveness of hydrodynamics on their migrations.

Oceanodromous fishes such as *Clupea harengus* (Herring), *Gadus morhua* (Cod), *Germo alalunga* (white tuna) and *Thunnus thynnus* (Atlantic Bluefin Tuna) live and migrate throughout the sea; Anadromous fishes such as *Salmo, Oncorhynchus* (salmon), live in the sea and migrate to cold, clear water of lakes or upstream rivers' gravel beds to breed; Catadromous fishes such as North American eel and European eel spend most of their lives in fresh water, then migrate to the sea to breed (Dorst 2019). On the other hand, potamodromous fishes such as salmonids and sturgeons shape and link among food webs, they are central aquatic species of environment and ecosystem of fresh water systems (Yoshioka 2017). As many of these fishes are economically valuable, the abundance and scarcity would cause a significant economic impact. Furthermore, it also affects

✉ Paramahansa Pramanik
  ppramanik1@niu.edu

1   Department of Mathematical Sciences, Northern Illinois University, 1425 Lincoln Highway, DeKalb, IL, USA

their habitats' ecosystems and the route of migration (Guse et al. 2015; Radinger and Wolter 2015; Yoshioka 2017).

A model describing the behavioral trade-off between migration time and energy expenditure has been discussed in Pinti et al. (2020). This model identifies optimal migration routes in realistic ocean conditions, and it explicitly includes a behavioral factor for individual risk management, including risks associated with moving in a stochastic oceanic environment, and it has been found that behavioral traits have significant influence in determining optimal routes in long-distance sea turtle migrations (Pinti et al. 2020). An extensive review of major models for animal migration such as analytic models, gametheoretic models, stochastic dynamic programming models and individual-based models have been discussed in Bauer and Klaassen (2013). A Hamiltonian–Jacobi–Bellman quasi-variational inequality (HJBQVI) equation has been used to determine an optimal migration strategy to give a maximized minimal profit (Yoshioka 2019). Furthermore, based on theoretical aspects Yoshioka (2019) suggests that sub-additivity of the performance index critically affects the resulting strategy. This model is a useful tool for comprehension of animal migration under different biological and environmental conditions and, from the viewpoint of the stochastic impulse control, violation of the sub-additivity was indicated to be an essential element for non-trivial migration strategies where not all the population migrates at once (Yoshioka 2019). In Yoshioka and Yaegashi (2018) a mathematical model for the onset of fish migration has been introduced in the context of a stochastic optimal stopping theory. Their analysis results provide the conditions for residency and migration. Furthermore, numerical computation in this paper turns out to be computationally feasible (Yoshioka and Yaegashi 2018).

Following Kappen (2007), every animal knows how to breath, digest, do elementary process of sensory information and motor action by birth. Therefore, it is automaton in nature. Animals are those type of automatons who learn from their environments and gain experiences about certain events such that in future they react in more intelligent manner. When a fish migrates to the spawning ground, its movements are influenced by unprecedented environmental noises. Furthermore, the migration processes are often exposed to influences that are incompletely understood. Therefore, extending deterministic or models incorporated with ordinary differential equations to ones that embrace more complicated variations are needed (Ton et al. 2014). From this viewpoint, one can introduce stochastic influences or noises. In general, as population dynamics in natural environment are always stochastic in nature (Lande et al. 2003), this leads the decision-making processes of the population under uncertainties (Yoshioka 2019). Therefore, stochastic models are appropriate in describing the population dynamics (Yoshioka 2019). Consider a fish is moving to the spawning ground with a school. On the way it faces a severe storm which makes very hard for that fish to stay in the school. Similarly, for an upstream migration if a fish got attacked by a predator, it is impossible for that fish to keep up with the velocity of its school. These environmental noises are uncertain in nature, and the fish's decision to stay in the school depends on it. The movement of a fish school is ergodic in nature, and therefore, its movement with some positive velocity is considered as a movement of a particle on a surface (Yoshioka 2017). Finally, as a fish is a very small part of a fish school, its movement is considered as the movement of a quantum particle in the same surface which I later describe as $\sqrt{8/3}$ Liouville-like quantum gravity surface. Based on different environmental noises, a fish determines its optimal strategy which leads to the shortest path out of infinite paths from the initial position to the terminal position by a Feynman-type path integral method (Feynman 1948).

Figure 1 gives some of all possible trajectories from the initial state $x_0$ at time $s = 0$ (i.e., initial position of the fish before migration). In this figure each trajectory represents the path of the fish from $x_0$ to the spawning area which it reaches at time $t$. The area inside the two arc-shaped lines is the feasible set of the optimal trajectory of the fish. At some point of time, if one trajectory goes beyond the feasible set then by Lebesgue–Riemann lemma it will come back inside in later time such that the path integration becomes measurable. Blue trajectory line goes beyond the feasible set just after starting at $x_0$. Red line is the optimal trajectory $\{x^o(s)\}$ for state $x$ (position of a fish) which ends up at terminal state $x(t)$ at time $t$. Almost every trajectory is upward sloping; therefore, action function can be used as a power of an exponential function. As forward looking method are considered here, there is no terminal condition and at each time point $s$ the fish has the information up to $s$.

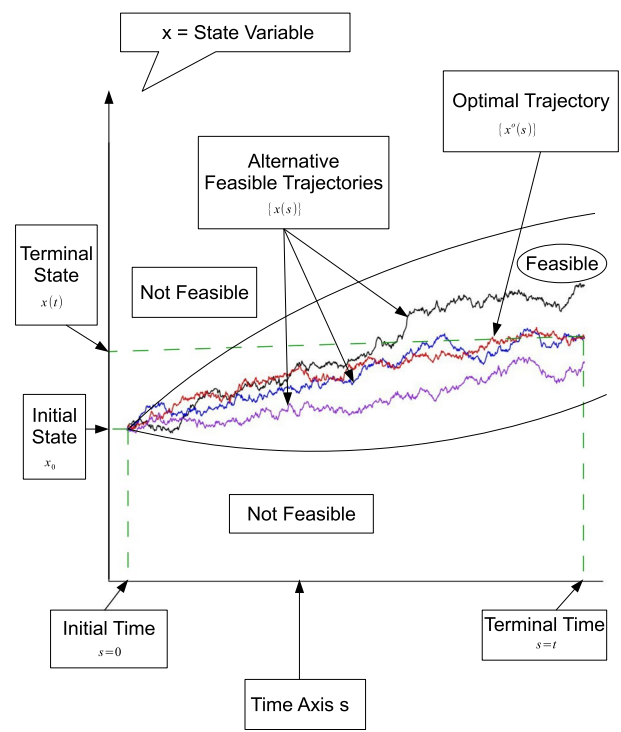

**Fig. 1** All possible trajectories of the fish during migration

When a fish swims against water current for spawning or foraging for food it faces unprecedented obstacles such as predators, change in climates, other obstacles in terms of change in path of a river or human made dams. If a fish did not face these obstacles before, it does not have enough information to react to it. On the other hand, if a fish survives the obstacle, then in future it knows how to react. In this paper three main behaviors of a fish school have been assumed; there is no leader in the fish school and each fish behaves same, each fish uses some form of weighted average of position and orientation of its nearest neighbors in order to decide its movement and, on the way to the destination it gains imperfect information with some degree of uncertainty which reflects to their actions (Nguyen et al. 2016). Furthermore, as the water current changes its direction and power instantaneously due to environmental factors, I further assume that the fish uses weighted average on a $\sqrt{8/3}$ Liouville-like quantum gravity surface (LQG) which can be glued to a Brownian surface of action of that fish (Duplantier and Sheffield 2011). The main reason of this assume is that each fish has its own action space with a dynamic strategy polygon which changes its shape instantaneously based on the available information. Furthermore, it creates a curve around itself on its action space. For example, if a predator is very near to the fish, it has to create an escape strategy such that the predator falls for the curvature in fish's action space with some probability, otherwise the fish dies. Following Nguyen et al. (2016) four behaviors of a fish school due to obstacles are assumed such as rebound, pullback, pass and reunion and separation. Furthermore, as each fish is assumed to be a quantum particle, these behaviors of a fish school are on $\sqrt{8/3}$-LQG surface.

The most important biological learning is Hebbian learning which states that, if two neurons become active simultaneously then the synaptic connections are stronger between them and slowly become weaker otherwise (Hebb 1949; Kappen 2007). The important assumption of this learning is that two neurons are uncorrelated to each other. There is a strong evidence that Hebbian learning occurs at *hippocampus* in the brain but it is too simple to consider in general behavior of synaptic plasticity (Kappen 2007). Many tasks are more complex and complicated than simultaneous Hebbian learning. They require imperfect information, uncertainty with some degrees, some sequential responses based on previous experiences and the result is only known at the future time. Examples of these tasks are motor control, foraging for food (Kappen 2007), and furthermore, for fish schools finding a spawning area.

A typical example for a motor control would be a fish trying to survive from a predator. Consider a fish is moving upstream and there is a bear at some point of the stream. At time $s$ the fish has to overcome this obstacle. Assume the location of the fish at time $s$ as the initial condition and reaching to the bear is the terminal condition. Then, the fish has infinite paths to come to the bear and it is a success for that bear in terms of catching fish. As both the bear and fish are automatons, when the bear catches the fish, it decodes another machine. As there are infinitely many paths between the initial and terminal conditions, to get an optimal path Feynman-type path integral control can be used (Feynman 1948). For any small $\epsilon > 0$, at time $s + \epsilon$ a fish's objective is not to come closer to the bear. Therefore, its terminal conditions is a horizon without the points where location of the bear is detected. As the bear is also moving over time to get a fish efficiently, the fish has to consider this fact while calculating its terminal conditions. A motor program is a sequence of actions: a path cost which specifies the energy consumption to contract the muscles in order to execute the motor program, and an end cost specifies whether the fish would come closer to the bear and get killed, just get hurt and escape, or manage to escape completely. Therefore, an optimal control solution is a sequence of motor commands that results in escaping of a fish from a bear which depends on the state and explicit on time. Similar situation happens in the ocean where a fish swimming against the ocean current get attacked by a great white shark. For detailed discussion see Kappen (2007).

When a fish is foraging for food or a spawning area, it explores the environment with the objective is to find as much food as possible in a short time period or, after it leaves another fish school might come and finish all the food. Similarly, if a first fish school does not exploit the area then another fish school might make it as a breeding ground. As availability of these places are purely random, at each time $s$, a fish considers the food it expects to encounter or find a place to breed in the period $[s, s + \epsilon]$. Now time horizon for a fish recedes into the future with the current time and the cost contributes a path with no end-cost (Kappen 2007). Hence, in each time point the fish faces same task at different locations of the environment which makes optimal control of that fish time-independent.

As each fish in a school is assumed to be a quantum particle I introduce an alternative method based on Feynman-type path integral to solve this stochastic control problem based on Feynman-type path integrals instead of traditional Pontryagin maximum principle. If the objective function is quadratic and the differential equations are linear, then solution is given in terms of a number of Ricatti equations which can be solved efficiently (Kappen 2007). But the water hydrodynamics is more complicated than just an ordinary linear differential equation and nonlinear stochastic feature gives the optimal solution a weighted mixture of suboptimal solutions, unlikely in the cases of deterministic or linear optimal control where a unique global optimal solution exists (Kappen 2007). In the presence of Wiener noise, Pontryagin maximum principle, a variational principle that leads to a coupled system of stochastic differential equations with initial and terminal conditions gives a generalized solution (Kappen 2007; Øksendal and Sulem 2019). Although incorporate randomness with its Hamiltonian–Jacobi–Bellman (HJB) equation is straight forward but difficulties come due to dimensionality when a numerical solution is calculated for both of deterministic or stochastic HJB (Kappen 2007). General stochastic control problem is intractable to solve computationally as it requires an exponential amount of memory and computational time because, the state space needs to be discretized and hence, becomes exponentially large in the number of dimensions (Theodorou et al. 2010; Theodorou 2011; Yang et al. 2014). Therefore, in order to calculate the expected values it is necessary to visit all states which leads to the summations of exponentially large sums (Kappen 2007; Yang et al. 2014). Kappen (2005a) and Kappen (2005b) say that a class of continuous nonlinear stochastic finite time horizon control problems can be solved more efficiently than Pontryagin maximum principle. These control problems reduce to computation of path integrals interpreted as free energy because, of their various statistical mechanics forms such as Laplace approximations, Monte Carlo sampling, mean field approximations or belief propagation (Kappen 2005a, b, 2007; Van Den Broek et al. 2008). According to Kappen (2007) these approximate computations are really fast.

Furthermore, a class of nonlinear HJB equations can be transformed into linear equations by doing a logarithmic transformation. This transformation stems back to the early days of quantum mechanics and was first used by Schrödinger to relate HJB equation to the Schrödinger equation (Kappen 2007). Because of this linear feature, backward integration of HJB equation over time can be replaced by computing expectation values under a forward diffusion process which requires a stochastic integration over trajectories that can be described by a path integral (Kappen 2007). Furthermore, in more generalized case like Merton–Garman–Hamiltonian system, getting a solution through Pontryagin maximum principle is impossible and Feynman path integral method gives a solution (Baaquie 1997). Previous works using the Feynman path integral method has been done in motor control theory by Kappen (2005b), Theodorou et al. (2010) and Theodorou (2011). A rigorous discussion of this quantum approach in finance has been done in Baaquie (2007). In Pramanik (2020) a Feynman-type path integral has been introduced to determine a feedback stochastic control. This method works in both linear and nonlinear stochastic differential equations and a Fourier transformation has been used to find out solution of Wick-rotated Schrödinger type equation (Pramanik 2020). This approach with $\sqrt{8/3}$-LQG is the first attempt to obtain an optimal strategy in the literature of fish migration.

## The model

Life cycle of a fish consists of growth, migration and reproduction under no overlapping generation (Yoshioka 2019). Furthermore, in Yoshioka (2019) optimal migration strategy and the basic reproduction number of an amphidromous fish *Plecoglossus altivelis altivelis* in Japan has been considered because, it is one of the most ecologically and commercially important fresh water fish. Apart from that, this fish has one-year life history seasonally migrating between sea and river (Yoshioka 2019). They grow up in the river during spring to coming summer by feeding algae like diatoms and go down-stream of the river to the sea for spawning. After the death of an adult fish, hatched larvae move toward the sea and grow up by feeding plankton until the next spring when usually the mass migration takes place (Yoshioka 2019). In this paper optimal migration strategy has been considered from one habitat to another by using a Feynman-type path integral approach, which is useful for nonlinear dynamics. Fishes like *Salangichthys microdon* and *Hypomesus nipponensis* have the similar spawning behavior like *Plecoglossus altivelis altivelis* (Arai et al. 2003, 2006; Yoshioka 2019).

Consider a fish moves from habitat $H_0$ to $H_1$ in $[0, t]$ period of time such that, $t > 0$. In this paper partial migration of a fish school is not considered, and it is assumed that an adult fish is not coming back from $H_1$ to $H_0$. Furthermore, the shape of the school only depends on external factors but not because of the fish interactions in the school. For time $s$ dependent state variables $x^i(s)$ and $v^i(s)$, and the control variable $u^i(s)$, $i^{th}$ fish has the function $h_{01}^i[s, x^i(s), v^i(s), u^i(s)]$ with the initial condition $h_{01}^{i*} \geq 0$ for all $i \in \{1, 2, ..., I\} \in I'$ and $s \in [0, t]$. The objective function $h_{01}^i[s, x^i(s), v^i(s), u^i(s)]$ is twice differentiable with respect to time in order to satisfy Wick rotation, is continuously differentiable with respect to $i^{th}$ fish's strategy $u^i(s)$, non-decreasing in state variables $x^i(s)$ and $v^i(s)$, non-increasing in $u^i(s)$, and convex and continuous in all state variables and strategies (Mas-Colell et al. 1995; Pramanik and Polansky 2020b). $H_{01}^i(s) \in [0, 1]$ is defined as the survival of $i^{th}$ fish during migration from $H_0$ to $H_1$. Instead of taking zero and unity $H_{01}^i(s)$ takes values in between them, because $i^{th}$ fish might be attacked by a predator and get severely injured, and that severity can be determined some number in between 0 and 1. Finally, at initial time 0 fish $i$ does not have any future information over $[0, t]$, it only makes expectations conditioned on initial states $x_0^i$ and $v_0^i$ at time 0 which is denoted as $\mathbb{E}_0$. Once the path integral method is introduced, entire time interval $[0, t]$ would be divided into smaller equaled time subintervals $[s, s + \epsilon]$ for all $\epsilon > 0$ such that $\epsilon \downarrow 0$ and $i^{th}$ fish makes expectation conditioned on its states $x^i(s)$ and $v^i(s)$ at time $s$ (hence, $\mathbb{E}_s$). The objective of fish $i \in \{1, 2, ..., I\} \in I'$ is :

$$
\begin{aligned}
\mathbf{OB}_\alpha^i &: \overline{\varPhi}_\alpha^i(s, x^i, v^i) \\
&= h_{01}^{i*} + \max_{u^i \in U} \mathbb{E}_0 \Bigg\{ \int_0^t \sum_{i=1}^I \exp(-\rho^i s) \alpha^i H_{01}^i(s) h_{01}^i \\
&\qquad [s, x^i(s), v^i(s), u^i(s)] \Big| \mathscr{F}_0^{x,v} \Bigg\} ds,
\end{aligned} \tag{1}
$$

where $u^i$ is the strategy of fish $i$ (control variable), $\alpha^i \in \mathbb{R}$ is constant weight, $I$ is total number of fishes in a school, $\rho_s^i \in (0, 1)$ is a stochastic discount rate of $i^{th}$ fish with $u^i \in \mathbb{R}^{I \times I}$, $v^i \in \mathbb{S}_i^{(I \times I) \times t}$ and $x^i \in \mathbb{R}^{I \times I}$ are time $s \geq 0$ dependent all possible controls and states available to them with two spheres available to $i^{th}$ fish is $\mathbb{S}_i^{(I \times I) \times t}$, and $\mathscr{F}_s^{x,v}$ is the $u^i$-adapted filtration process of hydrodynamics starting at the beginning of the migration process. For a fish the example of state variable $x^i$ might be its weight with more weight leads to more survivability during migration from $H_0$ to $H_1$ and become successful in spawning (Yoshioka 2019). Another state variable $v^i$ is its relative velocity against the water current. On the other hand, the control variable of a fish might be size of the school, strategy of a fish to get to $H_1$ and time spent at the spawning area (Yoshioka 2017, 2019). Without loss of generality, $x^i(s)$ and $v^i(s)$ are two state variables and $u^i(s)$ is assumed to be the control variable.

Now $i^{th}$ fish faces two stochastic hydrodynamic systems. The first system is

$$
\begin{aligned}
dx^i(s) = &\mu_1^i[s, x^i(s), v^i(s), u^i(s)]ds \\
&+ \sigma_1^i[s, x^i(s), v^i(s), u^i(s)]dB_1^i(s),
\end{aligned} \tag{2}
$$

where $\mu_1^i \in \mathbb{R}^{I \times I}$ is the drift component, $\sigma_1^i \in \mathbb{R}^{I \times I}$ is the diffusion component and $B_1^i(s)$ is an $I \times I$-dimensional standard Brownian motion. Similarly for the second-state variable $v^i(s)$ the dynamics is

$$
\begin{aligned}
dv^i(s) = &\mu_2^i[s, \psi, x^i(s), v^i(s), u^i(s)]ds \\
&+ \sigma_2^i[s, x^i(s), v^i(s), u^i(s)]dB_2^i(s),
\end{aligned} \tag{3}
$$

where $\mu_2^i$ is an $I \times I$-dimensional drift component, $\sigma_2^i$ is an $I \times I$-dimensional diffusion component, $B_2^i$ is the $I \times I$-dimensional standard Brownian motion process and the communication rate function between $i^{th}$ and $j^{th}$ fishes $\psi : [0, \infty) \to [0, \infty)$ is locally Lipschitz continuous (Nguyen et al. 2016). Clearly, in Eq. (2) if $\sigma_1^i = 0$, $\mu_1^i = v^i(s)$ and in Eq. (3) $\mu_2^i = \frac{\lambda}{I} \sum_{i=1}^I \psi(||x^i(s) - x^j(s)||)(v^j - v^i)$ with $\sigma_2^i = \sigma *$, where $\lambda$ is constant, nonnegative coupling strength

between two fishes (Ha et al. 2009), then the whole system is called Cucker–Smale system with white noise (Ahn and Ha 2010; Carrillo et al. 2010; Nguyen et al. 2016). Equations (2) and (3) represent more generalized version of Cucker–Smale system under white noise.

## Definitions and assumptions

**Assumption 1** For $t > 0$, suppose fish $i$ has drift components $\mu_1^i(s, x^i, v^i, u^i) : [0, t] \times \mathbb{R}^{I \times I} \times \mathbb{S}_i^{(I \times I) \times t} \times \mathbb{R}^{I \times I} \to \mathbb{R}^{I \times I}$, $\mu_2^i(s, \psi, x^i, v^i, u^i) : [0, t] \times \mathbb{R}^I \times \mathbb{R}^{I \times I} \times \mathbb{S}_i^{(I \times I) \times t} \times \mathbb{R}^{I \times I} \to \mathbb{R}^{I \times I}$ and the diffusion components $\sigma_1^i(s, x^i, v^i, u^i) : [0, t] \times \mathbb{R}^{I \times I} \times \mathbb{S}_i^{(I \times I) \times t} \times \mathbb{R}^{I \times I} \to \mathbb{R}^{I \times I}$, $\sigma_2^i(s, x^i, v^i, u^i) : [0, t] \times \mathbb{R}^{I \times I} \times \mathbb{S}_i^{(I \times I) \times t} \times \mathbb{R}^{I \times I} \to \mathbb{R}^{I \times I}$ are measurable functions with $(I \times I) \times t$-dimensional two-sphere $\mathbb{S}_i^{(I \times I) \times t}$ and, for some positive constants $K_1^i$ and $K_2^i$, $u^i \in \mathbb{R}^{I \times I}$, $x^i \in \mathbb{R}^{I \times I}$, and $v^i \in \mathbb{S}_i^{(I \times I) \times t}$ we have linear growth as

$$
\begin{aligned}
&|\mu_1^i(s, x^i, v^i, u^i)| + |\sigma_1^i(s, x^i, v^i, u^i)| \\
&\qquad \leq K_1^i(1 + |x^i| + |v^i|), \\
&|\mu_2^i(s, \psi, x^i, v^i, u^i)| + |\sigma_2^i(s, x^i, v^i, u^i)| \\
&\qquad \leq K_2^i(1 + |x^i| + |v^i|),
\end{aligned}
$$

such that, there exists another positive, finite, constants $K_3^i$ and $K_4^i$, and for different state variables $\tilde{x}^i_{(I \times I) \times 1}$ and $\tilde{v}^i_{(I \times I) \times 1}$ such that the Lipschitz conditions,

$$
\begin{aligned}
&|\mu_1^i(s, x^i, v^i, u^i) - \mu_1^i(s, \tilde{x}^i, v^i, u^i)| \\
&\qquad + |\sigma_1^i(s, x^i, v^i, u^i) - \sigma_1^i(s, \tilde{x}^i, v^i, u^i)| \\
&\qquad \leq K_3^i |x^i - \tilde{x}^i|, \\
&|\mu_2^i(s, \psi, x^i, v^i, u^i) - \mu_2^i(s, \psi, x^i, \tilde{v}^i, u^i)| \\
&\qquad + |\sigma_2^i(s, x^i, v^i, u^i) - \sigma_2^i(s, x^i, \tilde{v}^i, u^i)| \\
&\qquad \leq K_4^i |v^i - \tilde{v}^i|,
\end{aligned}
$$

are satisfied and for $\tilde{x}^i \in \mathbb{R}^{I \times I}$ and $\tilde{v}^i \in \mathbb{R}^{I \times I}$

$$
\begin{aligned}
&|\mu_1^i(s, x^i, v^i, u^i)|^2 + \|\sigma_1^i(s, x^i, v^i, u^i)\|^2 \\
&\qquad \leq (K_3^i)^2(1 + |\tilde{x}^i|^2 + |\tilde{v}^i|^2), \\
&|\mu_2^i(s, \psi, x^i, v^i, u^i)|^2 + \|\sigma_2^i(s, x^i, v^i, u^i)\|^2 \\
&\qquad \leq (K_4^i)^2(1 + |\tilde{x}^i|^2 + |\tilde{v}^i|^2),
\end{aligned}
$$

where $\|\sigma_m^i(s, x^i, v^i, u^i)\|^2 = \sum_{k=1}^I \sum_{l=1}^I |\sigma_m^{kl}(s, x^i, v^i, u^i)|^2$ for all $m = 1, 2$.

In Assumption 1 the state variable relative velocity $v^i \in \mathbb{S}_i^{(I \times I) \times t}$ is assumed to be on a two-sphere $\mathbb{S}_i$ such that it is homeomorphic but not diffeomorphic and, hence, it is a Brownian sphere. In general, if the velocity of water faced by $i^{th}$ fish is $w^i$ and that fish's velocity in stagnant water is $v_s^i$, then the relative velocity $v^i$ is $v_s^i - w^i$, where $w^i$ is constant (Yoshioka [2017]). This is a strong assumption. Water current changes in terms of direction as well as velocity. These changes might occur due to sudden environmental events such as tornadoes, flash floods, landslides, earth quakes, gravitation, volcanic eruptions under the sea and Thermohaline circulation. As all these events are random and fish $i$ does have uncertainty about it, its relative velocity is on this Brownian two-sphere. A Brownian surface of fish $i$ is a random Riemann surface parameterized by a domain on two-sphere whose Riemann metric tensor is $\exp\{\sqrt{8/3}k^i(l)\}dv^i \otimes dv^i$, where $k^i$ is some variant of the Gaussian free field (GFF) on some domain on this two-sphere, $l$ is some coming from two-sphere $\mathbb{S}_i^{(I \times I) \times t}$ and $dv^i \otimes dv^i$ is a Euclidean metric tensor (Gwynne and Miller [2016]). This is called $\sqrt{8/3}$-LQG surface. If this surface on quantum two-sphere has Schramm–Loewner–Evolution with parameter 6 ($SLE_6$) (Schramm [2000]), then under certain conditions the state variables show some upward and downward jumps (Miller [2018]) which considers the jump diffusion arises due to environmental conditions such as tornadoes, flash floods, landslides, earth quakes, gravitation, volcanic eruptions under the sea and Thermohaline circulation. Furthermore, $\sqrt{8/3}$-LQG surface glues to a Brownian surface (Gwynne and Miller [2016]; Sheffield [2007]; Sheffield et al. [2016]). In Section 5 this type of surface will be discussed.

**Assumption 2** Fish $i$ faces a probability space $(\Omega, \mathscr{F}_s^{,v}, \mathscr{P})$ with sample space $\Omega$, $u^i$-adaptive filtration at time $s$ of state variables $x^i$ and relative velocity $v^i$ as $\{\mathscr{F}_s^{,v}\} \subset \mathscr{F}_s$, a probability measure $\mathscr{P}$ and two $I \times I$-dimensional $\{\mathscr{F}_s\}$ Brownian motions $B_1^i$ and $B_2^i$ where the strategy of $i^{th}$ fish $u^i$ is an $\{\mathscr{F}_s^{,v}\}$ adapted process such that Assumption 1 holds, for the feedback control measure of fishes there exists a measurable function $h^i$ such that $h^i : [0, t] \times C([0, t]) : \mathbb{R}^{I \times I} \times \mathbb{S}_i^{(I \times I) \times t} \to u^i$ for which $u^i(s) = h^i[x^i(s, u^i), v^i(s, u^i)]$ such that Eqs. (2) and (3) have a strong unique solution (Ross [2008]).

**Assumption 3**

(i). $\mathscr{Z} \subset \mathbb{R}^{I \times I}$ such that fish $i$ cannot go beyond set $\mathscr{Z}_i \subset \mathscr{Z}$ because of its limitations of swimming against water current and different obstacles present in the water including the presence of a potential predator. This immediately implies set $\mathscr{Z}_i$ is different for different fishes. If the size of the fish is big, it would have lesser limitations and can swim more.

(ii). The function $h_0^i : [0, t] \times \mathbb{R}^{2I} \times \mathbb{S}_i^{(I \times I) \times t} \to \mathbb{R}^{I \times I}$. Therefore, all fishes in a school at the beginning of migration have the objective function $h_0 : [0, t] \times \mathbb{R}^{I \times I} \times \mathbb{S}_i^{(I \times I) \times t} \times \mathbb{R}^{I \times I} \to \mathbb{R}^{I \times I}$ such that

$h_0^i \subset h_0$ in functional spaces and both of them are concave which is equivalent to Slater condition (Marcet and Marimon [2019]). Possibility of partial migration of a school has been omitted in this paper.

(iii). There exists an $\epsilon > 0$ with $\epsilon \downarrow 0$ for all $(x^i, v^i, u^i)$ and $i = 1, 2, ..., I$ such that

$$\mathbb{E}_0 \left\{ \int_0^t \sum_{i=1}^{I} \exp(-\rho^i s)\alpha^i H_{01}^i(s)h_{01}^i \right.$$
$$\left. [s, x^i(s), v^i(s), u^i(s)] \middle| \mathscr{F}_0^{,v} \right\} ds \geq \epsilon.$$

The swimming path of $i^{th}$ fish during migration is continuous and it is mapped from an interval to a space of continuous functions with initial (the place where the migration begins) and terminal (i.e., spawning area, place where it finds food) points. Suppose, at time $s$, $g(s) : [p, q] \to \mathscr{C}$ represents a path of the migration of $i^{th}$ fish with initial and terminal points $g(p)$ and $g(q)$, respectively, such that, the line path integral is $\int_{\mathscr{C}} f(\gamma)ds = \int_p^q f(g(s))|g'(s)|ds$, where $g'(s)$ is derivative with respect to $s$. This paper concentrates on functional path integrals where the domain of the integral is the space of functions (Pramanik [2020]). Functional path integrals are very popular in probability theory and quantum mechanics. In Feynman ([1948]) theoretical physicist Richard Feynman introduced a new kind of functional path integral (Feynman path integral) and popularized it in quantum mechanics. Furthermore, mathematicians develop the measurability of this integral and in recent years it has become popular in probability theory (Fujiwara [2017]). In quantum mechanics, when a particle moves from one point to another, between those points it chooses the shortest path out of infinitely many paths such that some of them touch the edge of the universe. After introducing equal length small time interval $[s, s + \epsilon]$ with $\epsilon > 0$ such that $\epsilon \downarrow 0$ and using Riemann–Lebesgue lemma if at time $s$ one particle touches the end of the universe, then at a later time point it would come back and go to the opposite side of the previous direction to make the path integral a measurable function (Bochner et al. [1949]). Similarly, fish $i$ has infinitely many paths in between the initial migration point and the spawning area and, out of them, it chooses the optimal path given by the constraints explained in Eqs. (2) and (3). Furthermore, the advantage of Feynman approach is that it can be used in both in linear and nonlinear stochastic differential equation systems where constructing of an HJB equation is impossible (Baaquie [2007]). In this paper a Feynman-type path integral under $\sqrt{8/3}$-LQG has been introduced where each fish is assumed to be a quantum particle and there is no study so far on this type of approach in the fish migration literature.

**Definition 1** Suppose, $\mathscr{L}[s, y(s), \dot{y}(s)] = (1/2)m\dot{y}(s)^2 - V(y)$ be the classical Lagrangian function of a particle in generalized coordinate $y$ with mass $m$ where $(1/2) m\dot{y}^2$ and $V(y)$ are kinetic and potential energies, respectively. Therefore, the transition function of Feynman path integral corresponding to the classical action function $Z = \int_0^T \mathscr{L}(s, y(s), \dot{y}(s))ds$ is defined as $\Psi(y) = \int_{\mathbb{R}} \exp\{Z\}\mathscr{D}_Y$, where $\dot{y} = \partial y/\partial s$ and $\mathscr{D}_Y$ is an approximated Riemann measure which represents the positions of a particle at different time points $s$ (Pramanik 2020).

Here $i^{th}$ fish's objective is to maximize Eq. (1) subject to Eqs. (2) and (3). Following Definition 1 the quantum Lagrangian at time $s$ of $[s, s + \epsilon]$ is

$$
\begin{aligned}
\mathscr{L} = \mathbb{E}_s \Bigg\{ & \sum_{i=1}^{l} \exp(-\rho^i s)\alpha^i H_{01}^i(s) h_{01}^i[s, x^i(s), v^i(s), u^i(s)] \\
& + \lambda_1 \big[\Delta x^i(s) - \mu_1^i[s, x^i(s), v^i(s), u^i(s)]ds \\
& - \sigma_1^i[s, x^i(s), v^i(s), u^i(s)]dB_1^i(s)\big] \\
& + \lambda_2 \big[\Delta v^i(s) - \mu_2^i[s, \psi, x^i(s), v^i(s), u^i(s)]ds \\
& - \sigma_2^i[s, x^i(s), v^i(s), u^i(s)]dB_2^i(s)\big] \Bigg\},
\end{aligned}
\tag{4}
$$

where $\lambda_1$ and $\lambda_2$ are time-independent quantum Lagrangian multipliers. As at the beginning of the small time interval $[s, s + \epsilon]$, fish $i$ does not have any future information, it makes expectations based on its two state variables $x^i$ and $v^i$. For a penalization constant $L_\epsilon > 0$ and for time interval $[s, s + \epsilon]$ such that $\epsilon \downarrow 0$ define a transition function from $s$ to $s + \epsilon$ as

$$
\begin{aligned}
\Psi_{s,s+\epsilon}^i(x^i, v^i) = \frac{1}{L_\epsilon} \int_{\mathbb{R}^{l \times d}} \\
\exp[-\epsilon\mathscr{A}_{s,s+\epsilon}(x^i, v^i)]\Psi_s^i(x^i, v^i)dx^i(s) \times dv^i(s),
\end{aligned}
\tag{5}
$$

where $\Psi_s^i(x^i, v^i)$ is the value of the transition function at time $s$ with the initial condition $\Psi_0^i(x^i, v^i) = \Psi_0^i$ and the action function of fish $i$ is,

$$
\begin{aligned}
& \mathscr{A}_{s,s+\epsilon}(x^i, v^i) \\
& = \int_s^{s+\epsilon} \mathbb{E}_v \\
& \quad \Bigg\{ \sum_{i=1}^{l} \exp(-\rho^i s)\alpha^i H_{01}^i(v) h_{01}^i[s, x^i(v), v^i(v), u^i(v)]dv \\
& \quad + g^i[v + \Delta v, x^i(v) + \Delta x^i(v), v^i(v) + \Delta v^i(v)] \Bigg\},
\end{aligned}
$$

where $g^i[v + \Delta v, x^i(v) + \Delta x^i(v), v^i(v) + \Delta v^i(v)] \in C^2([0, t] \times \mathbb{R}^{l \times d} \times \mathbb{S}^{(l \times d) \times t})$ such that,

$$
\begin{aligned}
& g^i[v + \Delta v, x^i(v) + \Delta x^i(v), v^i(v) + \Delta v^i(v)] \\
& = \lambda_1 \big[\Delta x^i(v) - \mu_1^i[v, x^i(v), v^i(v), u^i(v)]dv \\
& - \sigma_1^i[v, x^i(v), v^i(v)), u^i(v)]dB_1^i(v)\big] \\
& + \lambda_2 \big[\Delta v^i(v) - \mu_2^i[v, \psi, x^i(v), v^i(v), u^i(v)]ds \\
& - \sigma_2^i[v, x^i(v), v^i(v), u^i(v)]dB_2^i(v)\big].
\end{aligned}
$$

Here the action function has the notation $\mathscr{A}_{s,s+\epsilon}(x^i, v^i)$ which means within $[s, s + \epsilon]$ the action of fish $i$ depends on the state variables $x^i$ and $v^i$, and furthermore, I assume this system has a feedback structure. Therefore, state variables also depend on the strategy of $i^{th}$ fish (i.e., $u^i$) as well as the rest of the school. Same argument goes to the transition function $\Psi_{s,s+\epsilon}^i(x^i, v^i)$.

**Definition 2** For fish $i$ optimal state variable $x^{i*}(s)$, relative velocity $v^{i*}(s)$ and its continuous optimal strategy $u^{i*}(s)$ constitute a dynamic stochastic Equilibrium such that for all $s \in [0, t]$ the conditional expectation of the objective function is

$$
\begin{aligned}
& \mathbb{E}_0 \Bigg[ \int_0^t \sum_{i=1}^{l} \exp(-\rho^i s)\alpha^i H_{01}^i(s) h_{01}^i \\
& \quad [s, x^{i*}(s), v^{i*}(s), u^{i*}(s)] \Big| \mathscr{F}_0^{x^*, v^*} \Bigg] ds \\
& \geq \mathbb{E}_0 \Bigg[ \int_0^t \sum_{i=1}^{l} \exp(-\rho^i s)\alpha^i H_{01}^i(s) h_{01}^i \\
& \quad [s, x^i(s), v^i(s), u^i(s)] \Big| \mathscr{F}_0^{x, v} \Bigg] ds,
\end{aligned}
$$

with the hydrodynamics explained in Eqs. (2) and (3), where $\mathscr{F}_0^{x^*, v^*}$ is the optimal filtration starting at time 0 such that, $\mathscr{F}_0^{x^*, v^*} \subset \mathscr{F}_0^{x, v}$.

## Link between HJB equation and path integral

Without loss of generality, in Eq. (1) assume $h_{01}^{i*} = 0$. Therefore, for a small time interval $[s, \tau]$ where $\tau = s + \epsilon$ for all $\epsilon \downarrow 0$ the objective function becomes,

$$
\begin{aligned}
& \overline{\Phi}_\alpha^i(s, x^i, v^i) \\
& = \max_{u^i \in U} \mathbb{E}_s \Bigg\{ \int_s^\tau \sum_{i=1}^{l} \exp(-\rho^i v)\alpha^i H_{01}^i(v) h_{01}^i \\
& \quad [v, x^i(v), v^i(v), u^i(v)] \Big| \mathscr{F}_v^{x, v} \Bigg\} dv.
\end{aligned}
$$

Consider $\alpha^i$ and $H^i_{01}$ are two constants and the function $h^i_{01}$ is quadratic with respect to the strategy such that $h^i_{01} = h^i_{01}(s, x^i, v^i) - \frac{1}{2}(u^i)^2$. After defining an arbitrary constant $R = \alpha^i H^i_{01} \exp(-\rho^i s)$ at time $s$ above objective function of fish $i$ becomes,

$$\overline{\Phi}^i_\alpha(s, x^i, v^i) = \max_{u^i \in U} \mathbb{E}_s \left\{ \int_s^\tau \sum_{i=1}^I W(v, x^i, v^i) \right.$$
$$\left. - \frac{R}{2}(u^i)^2 \middle| \mathscr{F}^{x,v}_v \right\} dv, \tag{6}$$

where $W(v, x^i, v^i) = \exp(-\rho^i v) \alpha^i H^i_{01} h^i_{01}(s, x^i, v^i)$, $h^i_{01}$ be an arbitrary function and $\mathbb{E}_s$ is the conditional expectation at time $s$ conditioned on state variables $x^i(s)$ and $v^i(s)$. For Eqs. (2) and (3) assume $\sigma^i_1$ and $\sigma^i_2$ are two constants and define

$$\mu^i_1(s, x^i, v^i, u^i) = \mu^i_1(s, x^i, v^i) + u^i,$$
$$\mu^i_2(s, \psi, x^i, v^i, u^i) = \mu^i_2(s, \psi, x^i, v^i) + u^i.$$

Therefore, $i^{th}$ fish's objective is to maximize Eq. (6) subject to Eqs. (2) and (3) such that above conditions on drift and diffusion hold. After setting $\tau_i = s + \epsilon$ a Taylor series expansion can be performed on $\overline{\Phi}^i_\alpha(\tau, x^i(\tau), v^i(\tau))$ around $s$ with first order with respect to $s$ and second order with respect to $x^i$ and $v^i$. By using Itô's lemma and following Baaquie (1997) we get,

$$\overline{\Phi}^i_\alpha[\tau, x^i(\tau), v^i(\tau)]$$
$$= \overline{\Phi}^i_\alpha(s, x^i, v^i) + \frac{\partial}{\partial s}\overline{\Phi}^i_\alpha(s, x^i, v^i) ds$$
$$+ \mu_1(s, x^i, v^i, u^i)\frac{\partial}{\partial x^i}\overline{\Phi}^i_\alpha(s, x^i, v^i) ds$$
$$+ \mu^i_2(s, \psi, x^i, v^i, u^i)\frac{\partial}{\partial v^i}\overline{\Phi}^i_\alpha(s, x^i, v^i) ds$$
$$+ \frac{1}{2}(\sigma^i_1)^2\frac{\partial^2}{\partial(x^i)^2}\overline{\Phi}^i_\alpha(s, x^i, v^i) ds$$
$$+ 2\rho(\sigma^i_1)^3\frac{\partial^2}{\partial x^i \partial v^i}\overline{\Phi}^i_\alpha(s, x^i, v^i) ds$$
$$+ (\sigma^i_2)^2\frac{\partial^2}{\partial(v^i)^2}\overline{\Phi}^i_\alpha(s, x^i, v^i) ds. ]$$

After using the conditions on drift and diffusion above equation becomes,

$$-\frac{\partial}{\partial s}\overline{\Phi}^i_\alpha(s, x^i, v^i) = \max_{u^i \in U}\left\{ \sum_{i=1}^I \left[ W(s, x^i, v^i) - \frac{R}{2}(u^i)^2 \right] \right.$$
$$+ \mu^i_1(s, x^i, v^i)\frac{\partial}{\partial x^i}\overline{\Phi}^i_\alpha(s, x^i, v^i)$$
$$+ u^i\frac{\partial}{\partial x^i}\overline{\Phi}^i_\alpha(s, x^i, v^i) + \mu^i_2(s, \psi, x^i, v^i)\frac{\partial}{\partial v^i}\overline{\Phi}^i_\alpha(s, x^i, v^i)$$
$$+ u^i\frac{\partial}{\partial v^i}\overline{\Phi}^i_\alpha(s, x^i, v^i)$$
$$+ \frac{1}{2}\left[ (\sigma^i_1)^2\frac{\partial^2}{\partial(x^i)^2}\overline{\Phi}^i_\alpha(s, x^i, v^i) \right.$$
$$\left. \left. + 2\rho(\sigma^i_1)^3\frac{\partial^2}{\partial x^i \partial v^i}\overline{\Phi}^i_\alpha(s, x^i, v^i) + (\sigma^i_2)^2\frac{\partial^2}{\partial(v^i)^2}\overline{\Phi}^i_\alpha(s, x^i, v^i) \right] \right\}. \tag{7}$$

After solving for the right hand side of Eq. (7) optimal strategy of fish $i$ is obtained as

$$u^{i*} = \frac{1}{R}\left[ \frac{\partial}{\partial x^i}\overline{\Phi}^i_\alpha(s, x^i, v^i) \right.$$
$$\left. + \frac{\partial}{\partial v^i}\overline{\Phi}^i_\alpha(s, x^i, v^i) \right]. \tag{8}$$

Using the result obtained in Eq. (8) and after plugging in to Eq. (7) yields,

$$-\frac{\partial}{\partial s}\overline{\Phi}^i_\alpha(s, x^i, v^i) = W(s, x^i, v^i)$$
$$- \frac{1}{2R}\left[ \frac{\partial}{\partial x^i}\overline{\Phi}^i_\alpha(s, x^i, v^i) + \frac{\partial}{\partial v^i}\overline{\Phi}^i_\alpha(s, x^i, v^i) \right]^2$$
$$+ \mu^i_1(s, x^i, v^i)\frac{\partial}{\partial x^i}\overline{\Phi}^i_\alpha(s, x^i, v^i)$$
$$+ \frac{1}{R}\left[ \frac{\partial}{\partial x^i}\overline{\Phi}^i_\alpha(s, x^i, v^i) \right]^2$$
$$+ \frac{2}{R}\left[ \frac{\partial}{\partial x^i}\overline{\Phi}^i_\alpha(s, x^i, v^i)\frac{\partial}{\partial v^i}\overline{\Phi}^i_\alpha(s, x^i, v^i) \right]$$
$$+ \mu^i_2(s, \psi, x^i, v^i)\frac{\partial}{\partial v^i}\overline{\Phi}^i_\alpha(s, x^i, v^i)$$
$$+ \frac{1}{R}\left[ \frac{\partial}{\partial v^i}\overline{\Phi}^i_\alpha(s, x^i, v^i) \right]^2$$
$$+ \frac{1}{2}\left[ (\sigma^i_1)^2\frac{\partial^2}{\partial(x^i)^2}\overline{\Phi}^i_\alpha(s, x^i, v^i) \right.$$
$$\left. + 2\rho(\sigma^i_1)^3\frac{\partial^2}{\partial x^i \partial v^i}\overline{\Phi}^i_\alpha(s, x^i, v^i) + (\sigma^i_2)^2\frac{\partial^2}{\partial(v^i)^2}\overline{\Phi}^i_\alpha(s, x^i, v^i) \right].$$

Stochastic HJB Eq. (9) is nonlinear with respect to $\overline{\Phi}^i_\alpha(s, x^i, v^i)$. The removal of the nonlinearity of Eq. (9)

would be a great help to solve this HJB Equation. After removing nonlinear parts HJB Eq. (9) becomes,

$$
\begin{aligned}
-\frac{\partial}{\partial s}\overline{\Phi}_\alpha^i(s,x^i,v^i) &= W(s,x^i,v^i) + \mu_1^i(s,x^i,v^i)\frac{\partial}{\partial x^i}\overline{\Phi}_\alpha^i(s,x^i,v^i) \\
&+ \frac{2}{R}\left[\frac{\partial}{\partial x^i}\overline{\Phi}_\alpha^i(s,x^i,v^i)\frac{\partial}{\partial v^i}\overline{\Phi}_\alpha^i(s,x^i,v^i)\right] \\
&+ \mu_2^i(s,\psi,x^i,v^i)\frac{\partial}{\partial v^i}\overline{\Phi}_\alpha^i(s,x^i,v^i) \\
&+ \frac{1}{2}\left[(\sigma_1^i)^2\frac{\partial^2}{\partial(x^i)^2}\overline{\Phi}_\alpha^i(s,x^i,v^i) + 2\rho(\sigma_1^i)^3\frac{\partial^2}{\partial x^i\partial v^i}\overline{\Phi}_\alpha^i(s,x^i,v^i)\right. \\
&+ \left.(\sigma_2^i)^2\frac{\partial^2}{\partial(v^i)^2}\overline{\Phi}_\alpha^i(s,x^i,v^i)\right].
\end{aligned}
\tag{10}
$$

D e f i n e, $\overline{\Phi}_\alpha^i(s,x^i,v^i) = -\omega\log\Theta(s,x^i,v^i)$, where $\omega = R[(\sigma_1^i)^2 + 2\rho(\sigma_1^i)^3 + (\sigma_2^i)^2]$ is nonzero. Hence,

$$
\begin{aligned}
\frac{\partial}{\partial s}\overline{\Phi}_\alpha^i(s,x^i,v^i) &= -\frac{\omega}{\Theta(s,x^i,v^i)}\left[\frac{\partial}{\partial s}\Theta(s,x^i,v^i,)\right], \\
\frac{\partial}{\partial x^i}\overline{\Phi}_\alpha^i(s,x^i,v^i) &= -\frac{\omega}{\Theta(s,x^i,v^i)}\left[\frac{\partial}{\partial x^i}\Theta(s,x^i,v^i)\right], \\
\frac{\partial}{\partial v^i}\overline{\Phi}_\alpha^i(s,x^i,v^i) &= -\frac{\omega}{\Theta(s,x^i,v^i)}\left[\frac{\partial}{\partial v^i}\Theta(s,x^i,v^i)\right], \\
\frac{\partial^2}{\partial(x^i)^2}\overline{\Phi}_\alpha^i(s,x^i,v^i) &\simeq -\frac{\omega}{\Theta(s,x^i,v^i)}\left[\frac{\partial^2}{\partial(x^i)^2}\Theta(s,x^i,v^i)\right], \\
\frac{\partial^2}{\partial(v^i)^2}\overline{\Phi}_\alpha^i(s,x^i,v^i) &\simeq -\frac{\omega}{\Theta(s,x^i,v^i)}\left[\frac{\partial^2}{\partial(v^i)^2}\Theta(s,x^i,v^i)\right], \\
\frac{\partial^2}{\partial x^i\partial v^i}\overline{\Phi}_\alpha^i(s,x^i,v^i) &\simeq -\frac{\omega}{\Theta(s,x^i,v^i)}\left[\frac{\partial^2}{\partial x^i\partial v^i}\Theta(s,x^i,v^i)\right].
\end{aligned}
\tag{11}
$$

Last three equations of Condition (11) are obtained by removing the quadratic part of $\Theta(s,x^i,v^i)$. Using Condition (11) yields,

$$
\begin{aligned}
-\frac{\partial}{\partial s}\Theta(s,x^i,v^i) &= -\frac{1}{\omega}W(s,x^i,v^i) \\
&+ \mu_1^i(s,x^i,v^i)\frac{\partial}{\partial x^i}\Theta(s,x^i,v^i) \\
&+ \mu_2^i(s,\psi,x^i,v^i)\frac{\partial}{\partial v^i}\Theta(s,x^i,v^i) \\
&+ \frac{1}{2}\left[(\sigma_1^i)^2\frac{\partial^2}{\partial(x^i)^2}\Theta(s,x^i,v^i)\right. \\
&+ 2\rho(\sigma_1^i)^3\frac{\partial^2}{\partial x^i\partial v^i}\Theta(s,x^i,v^i) \\
&+ \left.(\sigma_2^i)^2\frac{\partial^2}{\partial(v^i)^2}\Theta(s,x^i,v^i)\right].
\end{aligned}
\tag{12}
$$

For $\Theta(s,x^i,v^i) = \exp\{-\epsilon\mathscr{A}_{s,s+\epsilon}(x^i,v^i)\}$ with $\epsilon = 1/\omega$ stochastic HJB Eq. (12) must be solved backward in time. The removal of the nonlinearity of the HJB equation leads to reverse the direction of the computation in the following way. Consider fish $i$'s diffusion process $\Psi_s^i(\tilde{x}^i,\tilde{v}^i) = \Psi_s^i(\tau,\tilde{x}^i,\tilde{v}^i|s,x^i,v^i)$ for all $\tau > s$ represented by a Wick-rotated Schrödinger equation or a Fokker–Planck equation,

$$
\begin{aligned}
\frac{\partial}{\partial\tau}\Psi_s^i &= -\frac{W}{\omega}\Psi_s^i - \mu_1^i\frac{\partial}{\partial\tilde{x}^i}\Psi_s^i - \mu_2^i\frac{\partial}{\partial\tilde{v}^i}\Psi_s^i \\
&+ \frac{1}{2}\left[(\sigma_1^i)^2\frac{\partial^2}{\partial(\tilde{x}^i)^2}\Psi_s^i\right. \\
&+ 2\rho(\sigma_1^i)^3\frac{\partial^2}{\partial\tilde{x}^i\partial\tilde{v}^i}\Psi_s^i + (\sigma_2^i)^2\frac{\partial^2}{\partial(\tilde{v}^i)^2}\Psi_s^i\right],
\end{aligned}
\tag{13}
$$

with $\Psi_s^i(s,\tilde{x}^i,\tilde{v}^i|s,x^i,v^i) = \delta(\tilde{x}^i - x^i,\tilde{v}^i - v^i)$ being a Dirac delta function. Define

$$
B(s,x^i,v^i) = \int_{\mathbb{R}^{l\times l}}\Psi_s^i(s,\tilde{x}^i,\tilde{v}^i|s,x^i,v^i)\Theta(\tau,\tilde{x}^i,\tilde{v}^i)d\tilde{x}^i\times d\tilde{v}^i.
$$

Clearly $B(s,x^i,v^i)$ is independent of $\tau$ in both of the stochastic HJB Eq. (12) and the Fokker–Planck Eq. (13). Evaluating $B(s,x^i,v^i)$ for $\tau = s$ gives $B(s,x^i,v^i) = \Theta(s,x^i,v^i)$. Further evaluation $B(s,x^i,v^i)$ for $\tau = s_{\mu_1\mu_2}$ yields,

$$
\begin{aligned}
B(s,x^i,v^i) = \int_{\mathbb{R}^{l\times l}}\Psi_s^i(s_{\mu_1\mu_2}\tilde{x}^i,\tilde{v}^i|s,x^i,v^i)\Theta(s_{\mu_1\mu_2},x^i,v^i)d\tilde{x}^i \\
\times d\tilde{v}^i.
\end{aligned}
$$

Hence,

$$
\begin{aligned}
\Theta(s,x^i,v^i) = \int_{\mathbb{R}^{l\times l}}\exp\{-\epsilon\mathscr{A}_{s,s+\epsilon}(\tilde{x}^i,\tilde{v}^i)\}\Psi_s^i(s_{\mu_1\mu_2},\tilde{x}^i,\tilde{v}^i|s,x^i,v^i)d\tilde{x}^i \\
\times d\tilde{v}^i.
\end{aligned}
$$

Finally, after introducing the penalizing constant $L_\epsilon > 0$ above expression becomes,

$$
\begin{aligned}
\Psi_{s,s+\epsilon}(x^i,v^i) = \frac{1}{L_\epsilon}\int_{\mathbb{R}^{l\times l}}\exp\{-\epsilon\mathscr{A}_{s,s+\epsilon}(\tilde{x}^i,\tilde{v}^i)\}\Psi_s^i(\tilde{x}^i,\tilde{v}^i)d\tilde{x}^i \\
\times d\tilde{v}^i,
\end{aligned}
$$

which is the same expression of Feynman-type path integral as in Eq. (5). Therefore, Feynman-type path integral considers one class of stochastic HJB Equation.

## $\sqrt{8/3}$ Liouville quantum gravity surface

Consider the movement of a fish in the ocean. If it moves on a straight line, because of earth's spherical shape its path of movement is a curved line. Furthermore, as the water has currents, waves, environmental factors like volcanic eruptions underneath it, earthquakes, presence of predators, the fish sees the trajectories of paths have complicated shapes, curved in wild and random ways. In random geometry if the location of a fish is known, one can assign probabilities to the location of subsequent points. As here I subdivide the path of $i^{th}$ fish in very small parts with small equal length time interval $[s, s + \epsilon]$ such that $\epsilon \downarrow 0$, the random paths of fish $i$ resembles with a Brownian motion as it is the scaling limit of random walk. In a series of papers physicist Alexander Polyakov did explain properties of two-dimensional Brownian surface which is termed as Liouville quantum gravity (Polyakov 1981, 1987, 1996; Knizhnik et al. 1988; Pitici 2018). Later Scott Sheffield and Jason Miller mathematically prove that Brownian map (which takes the distance between two points on a random surface) and Liouville quantum gravity (which calculates the area) are fundamentally same (Miller and Sheffield 2016a, b; Pitici 2018).

As the relative velocity of fish $i$ is on an LQG surface, therefore, it is a random Riemann surface parameterized by a domain $\mathbb{D} \subset \mathbb{S}_i^{(I \times I) \times t}$ with Riemann metric tensor $e^{\gamma k^i(l)} \mathrm{d}v^i \otimes d\hat{v^i}$, where $\gamma \in (0, 2)$, $k^i$ is some variant of the Gaussian free field (GFF) on $\mathbb{D}$ (i.e., GFF and some harmonic function), $l$ is some number coming from two-sphere $\mathbb{S}_i^{(I \times I) \times t}$ and $\mathrm{d}v^i \otimes d\hat{v^i}$ is a Euclidean metric tensor (Gwynne and Miller 2016). This paper assumes $\gamma = \sqrt{8/3}$ because, it corresponds to a uniformly random planer map. Although traditional research suggests that forming a school improves navigation performances (Torney et al. 2015), hearing perception (Larsson 2012) and foraging efficiency (Wang et al. 2016), it leads to infinity large number of fish population in that school (which is absurd) (Yoshioka 2017). On the other hand, school formation could cause negative effects through passage efficiency (Lemasson et al. 2014), information transfer (Shang and Bouffanais 2014) and competition among fishes (Yoshioka 2017). Because of these effects fish $i$ has incomplete and imperfect information about its migration path trajectories, the action space is quantum in nature and $i^{th}$ fish's decision is a point on its dynamic convex strategy polygon of that quantum action space, and $k^i : \mathbb{S}_i^{(I \times I) \times t} \to \mathbb{R}^{I \times I}$ is a distribution such that $i^{th}$ fish's action can be represented by different trajectories.

Furthermore, for fish $i$, $\sqrt{8/3}$-LQG surface of its relative velocity is an equivalent class of action on two-sphere $(D, k^i)$ such that $D \subset \mathbb{S}_i^{(I \times I) \times t}$ is open and $k^i$ is a distribution function which is some variant of a GFF (Gwynne and Miller 2016). The pairs $(D, k^i)$ and $(\tilde{D}, \tilde{k^i})$ are equivalent if there exists a conformal map $\zeta : \tilde{D} \to D$ such that, $\tilde{k^i} = k^i \circ \zeta + Q \log |\zeta'|$, where $Q = 2/\gamma + \gamma/2 = \sqrt{3/2} + \sqrt{2/3}$ (Gwynne and Miller 2016). As the whole system is assumed to be a feedback system, the strategy space of all fishes in a school where the action is taken based on $v^i$ has the property like $\sqrt{8/3}$-LQG surface because, fish $i$ has radius $r^i$ around themselves such that, if another fish in the same school comes closer to compete, it would be able to handle. Furthermore, if $r^i \downarrow 0$, fish has a complete control over other fishes in the same school under the assumption that all the fishes are homogenous and have same level of skills. Therefore, the strategy space closer to the $i^{th}$ (i.e., $r^i = 0$) bends toward itself in such a way that the surface can be approximated to a surface on a two-sphere and furthermore, as the movement on this space is stochastic in nature, it behaves like a Brownian surface with its convex strategy polygon changes its shape at every time point based on the condition of water velocity. At $r^i = 0$ the surface hits essential singularity and fish $i$ has infinite power to control over other fishes.

**Definition 3** The knowledge space of $i^{th}$ fish $(\Omega, \mathscr{F}_s^{x,v}, \mathbb{S}_i, I')$ such that $i \in \{1, 2, ..., I\} \in I'$, each equivalent class with Riemann metric tensor $e^{\sqrt{8/3}k^i(l)} \mathrm{d}v^i \otimes d\hat{v^i}$ is finite, countably infinite or uncountable is defined as purely $\sqrt{8/3}$-LQG knowledge space which is purely quantum in nature (for detailed discussion about purely atomic knowledge see Hellman and Levy (2019)).

**Definition 4** Suppose, for $i^{th}$ fish $\mathbf{Y}_i(s, u^i)$ is a vector of two state variables $x^i$ and $v^i$ and is non-homogeneous Fellerian semigroup on time $s$ in $\sqrt{8/3}$-LQG surface $\mathbb{R}^{I \times I} \times \mathbb{S}_i^{(I \times I) \times t}$. The infinitesimal generator $A$ of $\mathbf{Y}_i(s, u^i)$ is defined by,

$$Ah_{01}^i(y) = \lim_{s \downarrow 0} \frac{\mathbb{E}_s[h_{01}^i(\mathbf{Y}_i(s, u^i))] - h_{01}^i(\mathbf{Y}_i(u^i))}{s},$$

for $\mathbf{Y}_i \in \mathbb{R}^{I \times I} \times \mathbb{S}_i^{(I \times I) \times t}$ where $h_{01}^i : [0, t] \times \mathbb{R}^{I \times I} \times \mathbb{S}_i^{(I \times I) \times t} \times \mathbb{R}^{I \times I} \to \mathbb{R}^{I \times I}$ is a $C_0^2\left(\mathbb{R}^{I \times I} \times \mathbb{S}_i^{(I \times I) \times t}\right)$ function, $\mathbf{Y}_i$ has a compact support, and at $\mathbf{Y}_i(u^i) > \mathbf{0}$ the limit exists where $\mathbb{E}_s$ represents the soccer team's conditional expectation on state variable $x^i$ and relative velocity $v^i$ at time $s$. Furthermore, if the above Fellerian semigroup is homogeneous on times, then $Ah_{01}^i$ is the Laplace operator in this space (Pramanik and Polansky 2020a).

**Definition 5** For a Fellerian semigroup $\mathbf{Y}_i(s, u^i)$ for all $\epsilon > 0$, the time interval $[s, s + \epsilon]$ with $\epsilon \downarrow 0$, define a characteristic-like quantum operator on $\sqrt{8/3}$-LQG surface starting at time $s$ as

$$\mathscr{A}h_{01}^i(\mathbf{Y}_i) = \lim_{\epsilon \downarrow 0} \frac{\log \mathbb{E}_s[\epsilon^2 h_{01}^i(\mathbf{Y}_i(s, u^i))] - \log[\epsilon^2 h_{01}^i(\mathbf{Y}_i(u^i))]}{\log \mathbb{E}_s(\epsilon^2)},$$

for $\mathbf{Y}_i \in \mathbb{R}^{I \times I} \times \mathbb{S}_i^{(I \times I) \times t}$, where $h_{01}^i : [0, t] \times \mathbb{R}^{I \times I} \times \mathbb{S}_i^{(I \times I) \times t} \to \mathbb{R}^{I \times I}$ is a $C_0^2\left(\mathbb{R}^{I \times I} \times \mathbb{S}_i^{(I \times I) \times t}\right)$ function, $\mathbb{E}_s$ represents the conditional expectation on $x^i$ and $v^i$ at time $s$, for $\epsilon > 0$ and a fixed $h_{01}^i$ we have the sets of all open balls of the form $B_\epsilon(h_{01}^i)$ contained in $\mathscr{B}$ (set of all open balls) and as $\epsilon \downarrow 0$ then $\log \mathbb{E}_s(\epsilon^2) \to \infty$.

**Assumption 4** The dynamic conditional expected objective function of $i^{th}$ fish explained in Eq. (1) on state variable dynamics $\mathbf{Y}_i \in \{\mathbf{Y}_i^0, \mathbf{Y}_i^1, ..., \mathbf{Y}_i^t\}$ is a tuple $\left(s, \boldsymbol{\alpha}^i, \{\mathbf{OB}_\alpha^i(u^i)\}_{u \in U}\right)$ where (i). $U$ is a finite strategy space based on two state variable feedback system where fish $i$ can choose strategy $u_i$ and $\boldsymbol{\alpha}^i$ is all probabilities available from which it chooses $\alpha^i$.

(ii). At time $s$, for each strategy $u_i \in U$, $\mathbf{OB}_\alpha^i \in \mathbb{R}^{\mathbf{Y}_i}$ is constrained objective function of fish $i$ such that Definition 5 holds.

## Main result

The components of stochastic differential game under $\sqrt{8/3}$-LQG with a continuum of states $x^i$ and $v^i$ and finite strategies available to fish $i$ are following:

– Let $\{1, 2, ..., I\} \in I'$ be a non-empty finite set of fishes in a school, $\mathscr{F}_s^{x,v}$ be the strategy adaptive filtration of state variable $x^i$ and relative velocity $v^i$ at time $s$ with the sample space $\Omega$.
– Fish $i$ has a finite set of strategies at time $s$ such that $u^i \in U$ for all $i \in I'$.
– Fish $i$ has discount rate $\rho^i \in (0, 1)$ with the constant weight $\alpha^i \in \mathbb{R}$.
– The bounded objective function $\mathbf{OB}_\alpha^i$ expressed in Eq. (1) is Borel measurable. Furthermore, this feedback system has a two dynamics expressed in Eqs. (2) and (3).
– As the relative velocity of $i^{th}$ fish is on $\mathbb{S}_i^{(I \times I) \times t}$, the migration process must be on a $\sqrt{8/3}$-LQG surface with the Riemann metric tensor $e^{\sqrt{8/3}k^i(l)}\mathrm{d}v^i \otimes \widehat{\mathrm{d}v^i}$, where $k^i$ is some variant of GFF (i.e., GFF and some harmonic function) such that $k^i : \mathbb{S}_i^{(I \times I) \times t} \to \mathbb{R}^{I \times I}$.
– For all $\varepsilon > 0$ with $\varepsilon \downarrow 0$ there exists a transition function from time $s$ to $s + \epsilon$ expressed as $\Psi_{s,s+\epsilon}^i(x^i, v^i) : \Omega \times \mathbb{R}^{I \times I} \times \mathbb{S}_i^{(I \times I) \times t} \times \prod_i u^i \to \varDelta\left(\Omega \times \mathbb{R}^{I \times I} \times \mathbb{S}_i^{(I \times I) \times t}\right)$ which is Borel-measurable.

The migration process for spawning or foraging for food is played in continuous time feedback environment. If $\mathbf{Y}_i \in \mathscr{Y} \subset (\Omega \times \mathbb{R}^{I \times I} \times \mathbb{S}_i^{(I \times I) \times t})$ be $i^{th}$ fish's condition on state variables before migration and it chooses a strategy profile at time $s$, $u_s^i$ such that $u_s^i \in \prod_i u^i$, then for $\varepsilon > 0$, $\Psi_{s,s+\epsilon}^i(\mathbf{Y}_i, u_s^i)$ is the conditional probability distribution of the next stage of the migration process. Fish $i$'s stable strategy depends on its relative velocity, obstacles including the presence of predators and behavior of the other members of the school at time $s$. Therefore, we can say it is Borel measurable mapping associates with state variable $\mathbf{Y}_i \subset \Omega$ a probability distribution on the set $u^i$.

**Proposition 1** *Suppose for all $i \in \{1, 2, ..., I\}$ fish $i$ has objective to maximize $\mathbf{OB}_\alpha^i$ with respect to $u^i \in U$ subject to two dynamical systems expressed in Equations (2) and (3) on*

$\sqrt{8/3}$-*LQG surface such that Assumptions 1-4 hold. Define a $C^2$ function*

$$
\begin{aligned}
f^i(s, x^i, v^i, u^i) &= \sum_{i=1}^I \exp(-\rho^i s)\alpha^i H_{01}^i h_{01}^i(s, x^i, v^i, u^i) \\
&\quad + g^i(s, x^i, v^i) \\
&\quad + \frac{\partial}{\partial s}g^i(s, x^i, v^i) + \frac{\partial}{\partial x^i}g^i(s, x^i, v^i) \\
&\quad \mu_1^i(s, x^i, v^i, u^i) \\
&\quad + \frac{\partial}{\partial v^i}g^i(s, x^i, v^i)\mu_2^i(s, \psi, x^i, v^i, u^i) \\
&\quad + \frac{1}{2}\bigg[\sigma_1^{i2}(s, x^i, v^i, u^i)\frac{\partial^2}{\partial x^i \partial x^{i'}}g^i(s, x^i, v^i) \\
&\quad + 2\rho\sigma_1^{i3}(s, x^i, v^i, u^i)\frac{\partial^2}{\partial x^i \partial v^i}g^i(s, x^i, v^i) \\
&\quad + \sigma_2^{i2}(s, x^i, v^i, u^i)\frac{\partial^2}{\partial v^i \partial v^{i'}}g^i(s, x^i, v^i)\bigg],
\end{aligned}
$$

*such that $g^i(s, x^i, v^i) \in C^2([0, t] \times \mathbb{R}^{I \times I} \times \mathbb{S}^{(I \times I) \times t})$ with Itô process $\widehat{Y}_i = g^i(s, x^i, v^i)$ is a positive, non-decreasing function vanishing at infinity. An optimal strategy of $i^{th}$ fish is the functional solution of*

$$
-\frac{\partial f^i(s, x^i, v^i, u^i)}{\partial u^i}\Psi_s^{i\tau}(x^i, v^i) = 0, \tag{14}
$$

*where a stable solution of $\Psi_s^{i\tau}(x^i, v^i)$ represented as*

$$
\Psi_s^i(x^i, v^i) = \exp\left\{-sf^i(s, x^i, v^i, u^i)\right\}\Psi_0^i(x^i, v^i)
$$

*is the transition wave function at time $s$ and states $x^i$ and $v^i$ with initial condition $\Psi_0^i(x^i, v^i) > 0$.*

**Proof** From quantum Lagrangian function expressed in Eq. (4) and after introducing $\sqrt{8/3}$-LQG surface for relative velocity $v^i$, Liouville–Feynman type action function of fish $i$ of time interval $[0, t]$ is

$$
\begin{aligned}
&\mathscr{A}_{0,t}(x^i, v^i) \\
&= \int_0^t \mathbb{E}_s\bigg\{\sum_{i=1}^I \exp(-\rho^i s)\alpha^i H_{01}^i(s)h_{01}^i \\
&\quad [s, x^i(s), v^i(s), u^i(s)] \\
&\quad + \lambda_1[\Delta x^i(s) - \mu_1^i[s, x^i(s), v^i(s), u^i(s)]\mathrm{d}s \\
&\quad - \sigma_1^i[s, x^i(s), v^i(s), u^i(s)]\mathrm{d}B_1^i(s)] \\
&\quad + \lambda_2[\Delta v^i(s) - \mu_2^i[s, \psi, x^i(s), v^i(s), u^i(s)]\mathrm{d}s \\
&\quad - \sigma_2^i[s, x^i(s), v^i(s), u^i(s)]\mathrm{d}B_2^i(s)] \\
&\quad + \lambda_3 e^{\sqrt{8/3}k^i(l(s))}\mathrm{d}s\bigg\},
\end{aligned}
$$

where $\lambda_1$, $\lambda_2$ and $\lambda_3$ are quantum Lagrangian time independent, nonnegative multipliers and furthermore, by $\lambda_3$ one can determine the presence of $\sqrt{8/3}$-LQG surface. Define $s + \epsilon = \tau$ such that $\epsilon > 0$ with $\epsilon \downarrow 0$, and for $L_\epsilon > 0$ convergence of path integral in Fujiwara (2017) implies

$$\Psi_{s,s+\epsilon}^i(x^i, v^i) = \frac{1}{L_\epsilon} \int_{\mathbb{R}^{2(I \times I)}}$$
$$\exp[-\epsilon \mathscr{A}_{s,s+\epsilon}(x^i, v^i)] \Psi_s^i(x^i, v^i) \mathrm{d}x^i(s) \times \mathrm{d}v^i(s), \quad (15)$$

where $\Psi_s^i(x^i, v^i)$ is the value of the transition function at time $s$ with the initial condition $\Psi_0^i(x^i, v^i) = \Psi_0^i$. As the time interval $[0, t]$ has been subdivided into $[s, \tau]$ equal lengthed small time intervals, the Liouville–Feynman type action function in that interval is

$$\mathscr{A}_{s,\tau}(x^i, v^i) = \int_s^\tau \mathbb{E}_v \Bigg\{ \sum_{i=1}^I \exp(-\rho^i v) \alpha^i H_{01}^i(v) h_{01}^i$$
$$[v, x^i(v), v^i(v), u^i(v)] dv$$
$$+ \lambda_1 \big[\Delta x^i(v) - \mu_1^i[v, x^i(v), v^i(v), u^i(v)] dv$$
$$- \sigma_1^i[v, x^i(v), v^i(v)), u^i(v)] dB_1^i(v) \big]$$
$$+ \lambda_2 \big[\Delta v^i(v) - \mu_2^i[v, \psi, x^i(v), v^i(v), u^i(v)] ds$$
$$- \sigma_2^i[v, x^i(v), u^i(v)] dB_2^i(v) \big]$$
$$+ \lambda_3 e^{\sqrt{8/3} k^i(l(v))} dv \Bigg\}, \quad (16)$$

with initial conditions $x^i(0) = x_0^i$ and $v^i(0) = v_0^i$, where $\Delta x^i(v) = x^i(v + dv) - x^i(v)$ and $\Delta v^i(v) = v^i(v + dv) - v^i(v)$. Now, Fubini's theorem implies,

$$\mathscr{A}_{s,\tau}(x^i, v^i) = \mathbb{E}_s \Bigg\{ \int_s^\tau \sum_{i=1}^I \exp(-\rho^i v) \alpha^i H_{01}^i(v) h_{01}^i$$
$$[v, x^i(v), v^i(v), u^i(v)] dv$$
$$+ \lambda_1 \big[\Delta x^i(v) - \mu_1^i[v, x^i(v), v^i(v), u^i(v)] dv$$
$$- \sigma_1^i[v, x^i(v), v^i(v)), u^i(v)] dB_1^i(v) \big]$$
$$+ \lambda_2 \big[\Delta v^i(v) - \mu_2^i[v, \psi, x^i(v), v^i(v), u^i(v)] ds$$
$$- \sigma_2^i[v, x^i(v), u^i(v)] dB_2^i(v) \big]$$
$$+ \lambda_3 e^{\sqrt{8/3} k^i(l(v))} dv \Bigg\}. \quad (17)$$

As $x^i(v)$ and $v^i(v)$ are Itô processes, Theorem 4.1.2 of Øksendal (2003) implies that there exists a function $g^i[v, x^i(v), v^i(v)] \in C^2([0, t] \times \mathbb{R}^{I \times I} \times \mathbb{S}^{(I \times I) \times t})$ such that Assumptions 1–4 hold and $\widehat{Y}_i(v) = g^i[v, x^i(v), v^i(v)]$, where $\widehat{Y}_i(v)$ is an Itô process. Assume

$$g^i[v + \Delta v, x^i(v) + \Delta x^i(v), v^i(v) + \Delta v^i(v)]$$
$$= \lambda_1 \big[\Delta x^i(v) - \mu_1^i[v, x^i(v), v^i(v), u^i(v)] dv$$
$$- \sigma_1^i[v, x^i(v), v^i(v)), u^i(v)] dB_1^i(v) \big]$$
$$+ \lambda_2 \big[\Delta v^i(v) - \mu_2^i[v, \psi, x^i(v), v^i(v), u^i(v)] dv$$
$$- \sigma_2^i[v, x^i(v), u^i(v)] dB_2^i(v) \big]$$
$$+ \lambda_3 e^{\sqrt{8/3} k^i(l(v))} dv + o(1).$$

Equation (17) becomes,

$$\mathscr{A}_{s,\tau}(x^i, v^i) = \mathbb{E}_s \Bigg\{ \int_s^\tau \sum_{i=1}^I$$
$$\exp(-\rho^i v) \alpha^i H_{01}^i(v) h_{01}^i$$
$$[v, x^i(v), v^i(v), u^i(v)] dv \qquad (18)$$
$$+ g^i[v + \Delta v, x^i(v) + \Delta x^i(v), v^i(v) + \Delta v^i(v)] \Bigg\}.$$

It is important to note that $g^i$ is not a function of either quantum Lagrangian multipliers (i.e., $\lambda_1$, $\lambda_2$ and $\lambda_3$) or the $\sqrt{8/3}$-LQG surface because, it takes those variables are parameters before the immigration process starts at time $s$. After using Itô's lemma and following Baaquie (1997), Eq. (18) becomes,

$$\mathscr{A}_{s,\tau}(x^i, v^i) = \sum_{i=1}^I \exp(-\rho^i s) \alpha^i H_{01}^i(s) h_{01}^i$$
$$[s, x^i(s), v^i(s), u^i(s)] + g^i[s, x^i(s), v^i(s)]$$
$$+ \frac{\partial}{\partial s} g^i[s, x^i(s), v^i(s)]$$
$$+ \frac{\partial}{\partial x^i} g^i[s, x^i(s), v^i(s)] \mu_1^i[s, x^i(s), v^i(s), u^i(s)]$$
$$+ \frac{\partial}{\partial v^i} g^i[s, x^i(s), v^i(s)] \mu_2^i[s, \psi, x^i(s), v^i(s), u^i(s)]$$
$$+ \frac{1}{2} \Big[ \sigma_1^{i2}[s, x^i(s), v^i(s), u^i(s)] \frac{\partial^2}{\partial x^i \partial x^{i'}} g^i[s, x^i(s), v^i(s)]$$
$$+ 2\rho \sigma_1^{i3}[s, x^i(s), v^i(s), u^i(s)] \frac{\partial^2}{\partial x^i \partial v^i} g^i[s, x^i(s), v^i(s)]$$
$$+ \sigma_2^{i2}[s, x^i(s), v^i(s), u^i(s)] \frac{\partial^2}{\partial v^i \partial v^{i'}} g^i[s, x^i(s), v^i(s)] \Big] + o(1), \quad (19)$$

where $\rho^2 < 1$ is the correlation coefficient between $x^i(s)$ and $v^i(s)$, $\sigma_1^{i2} = (\sigma_1^i)^2$, $\sigma_2^{i2} = (\sigma_2^i)^2$, $\sigma_1^{i3} = (\sigma_1^i)^3$, $x^{i'}(s)$ and $v^{i'}(s)$ are the transposition of the state variables $x^i(s)$ and $v^i(s)$, respectively. In Eq. (19), I have used the fact that $[\Delta x^i(s)]^2 = [\Delta v^i(s)]^2 = \epsilon$, and $\mathbb{E}_s[\Delta B_1^i(s)] = \mathbb{E}_s[\Delta B_2^i(s)]$, as $\epsilon \downarrow 0$ with initial conditions $x_0^i$ and $v_0^i$. Using Eq. (15), the transition wave function in $[s, \tau]$ becomes,

$$\Psi_{s,\tau}^i(x^i, v^i) = \frac{1}{L_\epsilon} \int_{\mathbb{R}^{2(I \times I)}}$$

$$\exp\Bigg\{ -\epsilon \Bigg[ \sum_{i=1}^{I} \exp(-\rho^i s) \alpha^i H_{01}^i(s) h_{01}^i$$

$$[s, x^i(s), v^i(s), u^i(s)]$$

$$+ g^i[s, x^i(s), v^i(s)] + \frac{\partial}{\partial s} g^i[s, x^i(s), v^i(s)]$$

$$+ \frac{\partial}{\partial x^i} g^i[s, x^i(s), v^i(s)] \mu_1^i[s, x^i(s), v^i(s), u^i(s)]$$

$$+ \frac{\partial}{\partial v^i} g^i[s, x^i(s), v^i(s)] \mu_2^i[s, \psi, x^i(s), v^i(s), u^i(s)]$$

$$+ \frac{1}{2} \Big[ \sigma_1^{i2}[s, x^i(s), v^i(s), u^i(s)]$$

$$\frac{\partial^2}{\partial x^i \partial x^{i'}} g^i[s, x^i(s), v^i(s)]$$

$$+ 2\rho \sigma_1^{i3}[s, x^i(s), v^i(s), u^i(s)]$$

$$\frac{\partial^2}{\partial x^i \partial v^i} g^i[s, x^i(s), v^i(s)]$$

$$+ \sigma_2^{i2}[s, x^i(s), v^i(s), u^i(s)] \frac{\partial^2}{\partial v^i \partial v^{i'}} g^i[s, x^i(s), v^i(s)] \Big] \Bigg] \Bigg\}$$

$$\Psi_s^i(x^i, v^i) \mathrm{d}x^i(s) \times \mathrm{d}v^i(s) + o(\epsilon^{1/2}),$$

(20)

as $\epsilon \downarrow 0$. For $\epsilon \downarrow 0$ define a new transition function $\Psi_s^{i\tau}(x^i, v^i)$ centered around time $\tau$ such that it can do the Taylor series expansion of $\Psi_{s,\tau}^i(x^i, v^i)$ up to order 1 in Eq. (20). Therefore,

$$\Psi_s^{i\tau}(x^i, v^i) + \epsilon \frac{\partial \Psi_s^{i\tau}(x^i, v^i)}{\partial s} + o(\epsilon)$$

$$= \frac{1}{L_\epsilon} \int_{\mathbb{R}^{2(I \times I)}} \exp\Bigg\{ -\epsilon$$

$$\Bigg[ \sum_{i=1}^{I} \exp(-\rho^i s) \alpha^i H_{01}^i(s) h_{01}^i[s, x^i(s), v^i(s), u^i(s)]$$

$$+ g^i[s, x^i(s), v^i(s)] + \frac{\partial}{\partial s} g^i[s, x^i(s), v^i(s)]$$

$$+ \frac{\partial}{\partial x^i} g^i[s, x^i(s), v^i(s)] \mu_1^i[s, x^i(s), v^i(s), u^i(s)]$$

$$+ \frac{\partial}{\partial v^i} g^i[s, x^i(s), v^i(s)] \mu_2^i[s, \psi, x^i(s), v^i(s), u^i(s)]$$

$$+ \frac{1}{2} \Big[ \sigma_1^{i2}[s, x^i(s), v^i(s), u^i(s)] \frac{\partial^2}{\partial x^i \partial x^{i'}} g^i[s, x^i(s), v^i(s)]$$

$$+ 2\rho \sigma_1^{i3}[s, x^i(s), v^i(s), u^i(s)] \frac{\partial^2}{\partial x^i \partial v^i} g^i[s, x^i(s), v^i(s)]$$

$$+ \sigma_2^{i2}[s, x^i(s), v^i(s), u^i(s)] \frac{\partial^2}{\partial v^i \partial v^{i'}} g^i[s, x^i(s), v^i(s)] \Big] \Bigg] \Bigg\}$$

$$\Psi_s^i(x^i, v^i) \mathrm{d}x^i(s) \times \mathrm{d}v^i(s) + o(\epsilon^{1/2}),$$

as $\epsilon \downarrow 0$.

For fixed $s$ and $\tau$ suppose that $x^i(s) = x^i(\tau) + \xi_1$, and $v^i(s) = v^i(\tau) + \xi_2$. For positive numbers $\eta_1 < \infty$ and $\eta_2 < \infty$ assume that $|\xi_1| \le \sqrt{\frac{\eta_1 \epsilon}{x^i(s)}}$ and $|\xi_2| \le \sqrt{\frac{\eta_2 \epsilon}{v^i(s)}}$. Here, two state variables of fish $i$ with the upper bounds are $x^i(s) \le \eta_1 \epsilon / \xi_1^2$ and $v^i(s) \le \eta_2 \epsilon / \xi_2^2$, respectively. Furthermore, by Fröhlich's reconstruction theorem (Simon 1979; Pramanik and Polansky 2020b) and Assumptions 1-4 imply

$$\Psi_s^{i\tau}(x^i, v^i) + \epsilon \frac{\partial \Psi_s^{i\tau}(x^i, v^i)}{\partial s} + o(\epsilon)$$

$$= \frac{1}{L_\epsilon} \int_{\mathbb{R}^{2(I \times I)}} \Bigg[ \Psi_s^{i\tau}(x^i, v^i) + \xi_1 \frac{\partial \Psi_s^{i\tau}(x^i, v^i)}{\partial x^i}$$

$$+ \xi_2 \frac{\partial \Psi_s^{i\tau}(x^i, v^i)}{\partial v^i} + o(\epsilon) \Bigg]$$

$$\times \exp\Bigg\{ -\epsilon \Bigg[ \sum_{i=1}^{I} \exp(-\rho^i s) \alpha^i H_{01}^i(s) h_{01}^i[s, x^i(\tau)$$

$$+ \xi_1, v^i(\tau) + \xi_2, u^i(s)]$$

$$+ g^i[s, x^i(\tau) + \xi_1, v^i(\tau)$$

$$+ \xi_2] + \frac{\partial}{\partial s} g^i[s, x^i(\tau) + \xi_1, v^i(\tau) + \xi_2]$$

$$+ \frac{\partial}{\partial x^i} g^i[s, x^i(\tau) + \xi_1, v^i(\tau) + \xi_2] \mu_1^i[s, x^i(\tau)$$

$$+ \xi_1, v^i(\tau) + \xi_2, u^i(s)]$$

$$+ \frac{\partial}{\partial v^i} g^i[s, x^i(\tau) + \xi_1, v^i(\tau) + \xi_2] \mu_2^i[s, \psi, x^i(\tau)$$

$$+ \xi_1, v^i(\tau) + \xi_2, u^i(s)]$$

$$+ \frac{1}{2} \Big[ \sigma_1^{i2}[s, x^i(\tau) + \xi_1, v^i(\tau) + \xi_2, u^i(s)]$$

$$\frac{\partial^2}{\partial x^i \partial x^{i'}} g^i[s, x^i(\tau) + \xi_1, v^i(\tau) + \xi_2]$$

$$+ 2\rho \sigma_1^{i3}[s, x^i(\tau) + \xi_1, v^i(\tau)$$

$$+ \xi_2, u^i(s)] \frac{\partial^2}{\partial x^i \partial v^i} g^i[s, x^i(\tau) + \xi_1, v^i(\tau) + \xi_2]$$

$$+ \sigma_2^{i2}[s, x^i(\tau) + \xi_1, v^i(\tau) + \xi_2, u^i(s)] \frac{\partial^2}{\partial v^i \partial v^{i'}} g^i[s, x^i(\tau)$$

$$+ \xi_1, v^i(\tau) + \xi_2] \Big] \Bigg] \Bigg\}$$

$$\Psi_s^i(x^i, v^i) d\xi_1 \times d\xi_2 + o(\epsilon^{1/2}),$$

(21)

as $\epsilon \downarrow 0$. For all $i \in \{1, 2, ..., I\}$ define a function

$$f^i[s, x^i, v^i, u^i(s)] = \sum_{i=1}^{I} \exp(-\rho^i s)\alpha^i H_{01}^i(s) h_{01}^i$$

$$[s, x^i(s), v^i(s), u^i(s)]$$

$$+ g^i[s, x^i(s), v^i(s)] + \frac{\partial}{\partial s} g^i[s, x^i(s), v^i(s)]$$

$$+ \frac{\partial}{\partial x^i} g^i[s, x^i(s), v^i(s)] \mu_1^i[s, x^i(s), v^i(s), u^i(s)]$$

$$+ \frac{\partial}{\partial v^i} g^i[s, x^i(s), v^i(s)] \mu_2^i[s, \psi, x^i(s), v^i(s), u^i(s)]$$

$$+ \frac{1}{2}\Big[ \sigma_1^{i2}[s, x^i(s), v^i(s), u^i(s)] \frac{\partial^2}{\partial x^i \partial x^{i\prime}} g^i[s, x^i(s), v^i(s)]$$

$$+ 2\rho \sigma_1^{i3}[s, x^i(s), v^i(s), u^i(s)] \frac{\partial^2}{\partial x^i \partial v^{i\prime}} g^i[s, x^i(s), v^i(s)]$$

$$+ \sigma_2^{i2}[s, x^i(s), v^i(s), u^i(s)] \frac{\partial^2}{\partial v^i \partial v^{i\prime}} g^i[s, x^i(s), v^i(s)] \Big].$$

Then Eq. (21) becomes,

$$\Psi_s^{i\tau}(x^i, v^i) + \epsilon \frac{\partial \Psi_s^{i\tau}(x^i, v^i)}{\partial s} + o(\epsilon)$$

$$= \frac{1}{L_\epsilon} \Psi_s^{i\tau}(x^i, v^i) \int_{\mathbb{R}^{2(I \times I)}}$$

$$\exp\{-\epsilon f^i[s, \xi_1, \xi_2, u^i(s)]\} d\xi_1 d\xi_2$$

$$+ \frac{1}{L_\epsilon} \frac{\partial \Psi_s^{i\tau}(x^i, v^i)}{\partial x^i} \int_{\mathbb{R}^{2(I \times I)}} \xi_1$$

$$\exp\{-\epsilon f^i[s, \xi_1, \xi_2, u^i(s)]\} d\xi_1 d\xi_2$$

$$+ \frac{1}{L_\epsilon} \frac{\partial \Psi_s^{i\tau}(x^i, v^i)}{\partial v^i} \int_{\mathbb{R}^{2(I \times I)}} \xi_2$$

$$\exp\{-\epsilon f^i[s, \xi_1, \xi_2, u^i(s)]\} d\xi_1 d\xi_2 + o(\epsilon^{1/2})$$

Assume that $f^i[s, \xi_1, \xi_2, u^i(s)]$ is a $C^2$ function, then

$$f^i[s, \xi_1, \xi_2, u^i(s)] = f^i[s, x^i(\tau), v^i(\tau), u^i(s)]$$

$$+ [\xi_1 - x^i(\tau)] \frac{\partial}{\partial x^i} f^i[s, x^i(\tau), v^i(\tau), u^i(s)]$$

$$+ [\xi_2 - v^i(\tau)] \frac{\partial}{\partial v^i} f^i[s, x^i(\tau), v^i(\tau), u^i(s)]$$

$$+ \frac{1}{2}\Big[ [\xi_1 - x^i(\tau)]'[\xi_1 - x^i(\tau)]$$

$$\frac{\partial^2}{\partial x^i \partial x^{i\prime}} f^i[s, x^i(\tau), v^i(\tau), u^i(s)]$$

$$+ 2[\xi_1 - x^i(\tau)][\xi_2 - v^i(\tau)]$$

$$\frac{\partial^2}{\partial x^i \partial v^{i\prime}} f^i[s, x^i(\tau), v^i(\tau), u^i(s)]$$

$$+ [\xi_2 - v^i(\tau)]'[\xi_2 - v^i(\tau)]$$

$$\frac{\partial^2}{\partial v^i \partial v^{i\prime}} f^i[s, x^i(\tau), v^i(\tau), u^i(s)] \Big] + o(\epsilon),$$

as $\epsilon \downarrow 0$ and $\Delta u^i \downarrow 0$, where $[\xi_1 - x^i(\tau)]'$ and $[\xi_2 - v^i(\tau)]'$ is the transposition of $[\xi_1 - x^i(\tau)]$ and $[\xi_2 - v^i(\tau)]$, respectively.

Define $m_1^i = \xi_1 - x^i(\tau)$ and $m_2^i = \xi_2 - v^i(\tau)$ so that $d\xi_1 = dm_1^i$ and $d\xi_2 = dm_2^i$, respectively, so that

$$\int_{\mathbb{R}^{2(I \times I)}} \exp\{-\epsilon f^i[s, \xi_1, \xi_2, u^i(s)]\} d\xi_1 d\xi_2$$

$$= \int_{\mathbb{R}^{2(I \times I)}} \exp\Big\{ -\epsilon \Big[ f^i[s, x^i(\tau), v^i(\tau), u^i(s)]$$

$$+ m_1^i \frac{\partial}{\partial x^i} f^i[s, x^i(\tau), v^i(\tau), u^i(s)]$$

$$+ m_2^i \frac{\partial}{\partial v^i} f^i[s, x^i(\tau), v^i(\tau), u^i(s)] \qquad (22)$$

$$+ \frac{1}{2} m_1^{i\prime} m_1^i \frac{\partial^2}{\partial x^i \partial x^{i\prime}} f^i[s, x^i(\tau), v^i(\tau), u^i(s)]$$

$$+ m_1^i m_2^i \frac{\partial^2}{\partial x^i \partial v^{i\prime}} f^i[s, x^i(\tau), v^i(\tau), u^i(s)]$$

$$+ \frac{1}{2} m_2^{i\prime} m_2^i \frac{\partial^2}{\partial v^i \partial v^{i\prime}} f^i[s, x^i(\tau), v^i(\tau), u^i(s)] \Big] \Big\} dm_1^i dm_2^i,$$

where $m_1^{i\prime}$ and $m_2^{i\prime}$ are the transposition of $m_1^i$ and $m_2^i$, respectively.

Let

$$\Theta^i = \begin{bmatrix} \frac{1}{2} \frac{\partial^2}{\partial x^i \partial x^{i\prime}} f^i[s, x^i(\tau), v^i(\tau), u^i(s)] & \frac{1}{2} \frac{\partial^2}{\partial x^i \partial v^{i\prime}} f^i[s, x^i(\tau), v^i(\tau), u^i(s)] \\ \frac{1}{2} \frac{\partial^2}{\partial x^i \partial v^{i\prime}} f^i[s, x^i(\tau), v^i(\tau), u^i(s)] & \frac{1}{2} \frac{\partial^2}{\partial v^i \partial v^{i\prime}} f^i[s, x^i(\tau), v^i(\tau), u^i(s)] \end{bmatrix},$$

and

$$m^i = \begin{bmatrix} m_1^i \\ m_2^i \end{bmatrix},$$

and

$$-V_1^i = \begin{bmatrix} \frac{\partial}{\partial x^i} f^i[s, x^i(\tau), v^i(\tau), u^i(s)] \\ \frac{\partial}{\partial v^i} f^i[s, x^i(\tau), v^i(\tau), u^i(s)] \end{bmatrix},$$

where I assume that $\Theta^i$ is positive definite, then the integrand in Eq. (22) becomes a shifted Gaussian integral,

$$\int_{\mathbb{R}^{2(I \times I)}} \exp\Big\{ -\epsilon \big( f^i - V_1^{i\prime} m^i + m^{i\prime} \Theta^i m^i \big) \Big\} dm^i$$

$$= \exp\big(-\epsilon f^i\big) \int_{\mathbb{R}^{2(I \times I)}} \exp\Big\{ (\epsilon V_1^{i\prime}) m^i - m^{i\prime}(\epsilon \Theta^i) m^i \Big\} dm^i$$

$$= \frac{\pi}{\sqrt{\epsilon |\Theta^i|}} \exp\Big[ \frac{\epsilon}{4} V_1^{i\prime} (\Theta^i)^{-1} V_1^i - \epsilon f^i \Big],$$

where $V_1^{i\prime}$ and $m^{i\prime}$ are the transposition of vectors $V_1^i$ and $m^i$, respectively. Therefore,

$$\frac{1}{L_\epsilon} \Psi_s^{i\tau}(x^i, v^i) \int_{\mathbb{R}^{2(I \times I)}} \exp\Big\{ -\epsilon f^i[s, \xi_1, \xi_2, u^i(s)] \Big\} d\xi_1 d\xi_2$$

$$= \frac{1}{L_\epsilon} \Psi_s^{i\tau}(x^i, v^i) \frac{\pi}{\sqrt{\epsilon |\Theta^i|}} \exp\Big[ \frac{\epsilon}{4} V_1^i (\Theta^i)^{-1} V_1^i - \epsilon f^i \Big], \qquad (23)$$

such that inverse matrix $(\Theta^i)^{-1} > 0$ exists. Similarly,

$$\frac{1}{L_\epsilon} \frac{\partial \Psi_s^{i\tau}(x^i, v^i)}{\partial x^i} \int_{\mathbb{R}^{2(l \times l)}} \xi_1$$
$$\exp\left\{-\epsilon f^i[s, \xi_1, \xi_2, u^i(s)]\right\} d\xi_1 d\xi_2$$
$$= \frac{1}{L_\epsilon} \frac{\partial \Psi_s^{i\tau}(x^i, v^i)}{\partial x^i} \frac{\pi}{\sqrt{\epsilon|\Theta^i|}} \left(\frac{1}{2}(\Theta^i)^{-1} + x^i\right)$$
$$\exp\left[\frac{\epsilon}{4} V_1^{i\prime} (\Theta^i)^{-1} V_1^i - \epsilon f^i\right], \tag{24}$$

and

$$\frac{1}{L_\epsilon} \frac{\partial \Psi_s^{i\tau}(x^i, v^i)}{\partial v^i} \int_{\mathbb{R}^{2(l \times l)}} \xi_1$$
$$\exp\left\{-\epsilon f^i[s, \xi_1, \xi_2, u^i(s)]\right\} d\xi_1 d\xi_2$$
$$= \frac{1}{L_\epsilon} \frac{\partial \Psi_s^{i\tau}(x^i, v^i)}{\partial v^i} \frac{\pi}{\sqrt{\epsilon|\Theta^i|}} \left(\frac{1}{2}(\Theta^i)^{-1} + v^i\right)$$
$$\exp\left[\frac{\epsilon}{4} V_1^{i\prime} (\Theta^i)^{-1} V_1^i - \epsilon f^i\right]. \tag{25}$$

Equations (23), (24) and (25) imply that the Wick-rotated Schrödinger type equation is,

$$\Psi_s^{i\tau}(x^i, v^i) + \epsilon \frac{\partial \Psi_s^{i\tau}(x^i, v^i)}{\partial s} + o(\epsilon)$$
$$= \frac{1}{L_\epsilon} \frac{\pi}{\sqrt{\epsilon|\Theta^i|}} \exp\left[\frac{\epsilon}{4} V_1^{i\prime} (\Theta^i)^{-1} V_1^i - \epsilon f^i\right]$$
$$\left[\Psi_s^{i\tau}(x^i, v^i) + \left(\frac{1}{2}(\Theta^i)^{-1} + x^i\right) \frac{\partial \Psi_s^{i\tau}(x^i, v^i)}{\partial K}\right.$$
$$\left. + \left(\frac{1}{2}(\Theta^i)^{-1} + v^i\right) \frac{\partial \Psi_s^{i\tau}(x^i, v^i)}{\partial v^i}\right] + o(\epsilon^{1/2}),$$

as $\epsilon \downarrow 0$.

Assuming $L_\epsilon = \pi/\sqrt{\epsilon|\Theta^i|} > 0$,

$$\Psi_s^{i\tau}(x^i, v^i) + \epsilon \frac{\partial \Psi_s^{i\tau}(x^i, v^i)}{\partial s} + o(\epsilon)$$
$$= \left[1 + \epsilon\left(\frac{1}{4} V_1^{i\prime} (\Theta^i)^{-1} V_1^i - \epsilon f^i\right)\right]$$
$$\left[\Psi_s^{i\tau}(x^i, v^i) + \left(\frac{1}{2}(\Theta^i)^{-1} + x^i\right) \frac{\partial \Psi_s^{i\tau}(x^i, v^i)}{\partial x^i}\right.$$
$$\left. + \left(\frac{1}{2}(\Theta^i)^{-1} + v^i\right) \frac{\partial \Psi_s^{i\tau}(x^i, v^i)}{\partial v^i}\right] + o(\epsilon^{1/2}). \tag{26}$$

As $x^i(s) \leq \eta_1 \epsilon / \xi_1^2$, assume $|(\Theta^i)^{-1}| \leq 2\eta_1 \epsilon(1 - \xi_1^{-2})$ such that $|(2\Theta^i)^{-1} + x^i| \leq \eta_1 \epsilon$. For $v^i(s) \leq \eta_2 \epsilon / \xi_2^2$ I assume $|(\Theta^i)^{-1}| \leq 2\eta_2 \epsilon(1 - \xi_2^{-2})$ such that $|(2\Theta^i)^{-1} + v^i| \leq \eta_2 \epsilon$. Therefore, $|(\Theta^i)^{-1}| \leq 2\epsilon \min\left\{\eta_1(1 - \xi_1^{-2}), \eta_2(1 - \xi_2^{-2})\right\}$ such that, $|(2\Theta^i)^{-1} + x^i| \to 0$ and $|(2\Theta^i)^{-1} + v^i| \to 0$. Hence

$$\Psi_s^{i\tau}(x^i, v^i) + \epsilon \frac{\partial \Psi_s^{i\tau}(x^i, v^i)}{\partial s} + o(\epsilon)$$
$$= (1 - \epsilon f^i)\Psi_s^{i\tau}(x^i, v^i) + o(\epsilon^{1/2}).$$

Therefore, the Wick-rotated Schrödinger type equation is,

$$\frac{\partial \Psi_s^{i\tau}(x^i, v^i)}{\partial s}$$
$$= -f^i[s, \xi_1, \xi_2, u^i(s)]\Psi_s^{i\tau}(x^i, v^i).$$

Therefore, the solution of

$$-\frac{\partial f^i[s, \xi_1, \xi_2, u^i(s)]}{\partial u^i} \Psi_s^i(x^i, v^i) = 0, \tag{27}$$

is an optimal strategy of fish $i$. Furthermore, as $\xi_1 = x^i(s) - x^i(\tau)$ and $\xi_2 = v^i(s) - v^i(\tau)$, for $\epsilon \downarrow 0$, in Eq. (27) $\xi_1$ and $\xi_2$ can be replaced by $x^i$ and $v^i$, respectively. Following Pramanik (2020) a stable solution to the wave function is

$$\Psi_s^i(x^i, v^i) = \exp\left\{-sf^i(s, x^i, v^i, u)\right\} \Psi_0^i(x^i, v^i),$$

where $\Psi_0^i(x^i, v^i)$ is the initial condition of the migration process of fish $i$. As the transition function $\Psi_s^i(x^i, v^i)$ is the solution to Eq. (27), the result follows. □

**Example 1** Following Yoshioka (2019) suppose, $i^{th}$ fish's objective is to maximize

$$\mathbf{OB}_\alpha^i : \overline{\Phi}_\alpha^i(s, x^i, v^i)$$
$$= h_{01}^{i*} + \max_{u^i \in U} \mathbb{E}_0 \left\{ \int_0^t \sum_{i=1}^I \exp(-\rho^i s) \alpha^i H_{01}^{i*}(s) x^i(s) v^i(s) [u^i(s)]^2 \right.$$
$$\left. \left| \mathscr{P}_0^{x,v} \right\} ds, \tag{28}$$

where fish $i$ assumed to be completely survived the migration process from habitat $H_0$ to $H_1$ and $H_{01}^i$ takes a constant value $H_{01}^{i*}$ and $h_{01}^i(s, x^i, v^i, u^i) = x^i v^i (u^i)^2$. Here assume $x^i(s)$ is the position of $i^{th}$ fish at time $s$ and $v^i(s)$ is the relative velocity of it. The objective function expressed in Yoshioka (2019) has one state variable and the control variable has been used without any exponent unlike in Eq. (28), where two state variables are used and the control variable has the exponent of two. Fish $i$ assumed to have strategies of its two state variables which makes $u^i$ come as a squared term in Eq. (28). First constraint of the fish is

$$dx^i(s) = u^i(s)v^i(s)ds + \sigma_1^i dB_1^i(s), \tag{29}$$

where $\sigma_1^i dB_1^i(s)$ is the noise resulting from the imperfectness of information-gathering and action of fish $i$ (Uchitane et al. 2012; Nguyen et al. 2016). Without loss of generality the diffusion component $\sigma_1^i$ is assumed to be constant. Stochastic differential equation represented by Eq. (29) and the equation represented in Uchitane et al. (2012) and Nguyen et al. (2016) are similar. Only $u^i$ has been added with relative velocity $v^i$. In order to derive $dv^i$ a Cucker–Smale type of system under white noise has been introduced (Nguyen et al. 2016) where the communication rate between $i^{th}$ and $j^{th}$ fishes $\psi : [0, \infty) \to [0, \infty)$ is assumed to be constant. Hence, the second constraint fish $i$ faces is

$$dv^i(s) = \frac{\lambda}{I} \sum_{i=1}^{I} u^i(s) \big( ||x^i(s) - x^j(s)|| \big) [v^i(s) - v^j(s)] ds$$
$$+ \sqrt{\sigma_2^{i*}} dB_2^i(s), \tag{30}$$

where $\lambda$ is constant, nonnegative coupling strength between two fishes (Ha et al. 2009), $||x^i(s) - x^j(s)||$ is assumed to be a Euclidean norm and $\sigma_2^{i*}$ is constant diffusion component (Ton et al. 2014). As I assume the system is a feedback system, before calculating an optimal strategy $u^i$ fish $i$ knows $x^i$ and $v^i$ at time $s$. Therefore, in Eqs. (28), (29) and (30), $x^i(s) = x^i$, $v^i(s) = v^i$ and $u^i(s) = u^i$. Following Øksendal (2003) assume

$$g^i(s, x^i, v^i) = \exp \left\{ s\lambda_1 v^i + \frac{s\lambda\lambda_2}{I} \sum_{i=1}^{I} \psi(x^i - x^j)(v^i - v^j) \right.$$
$$\left. + \lambda_3 \sqrt{8/3} k^i(l) \right\}.$$

Then

$$\frac{\partial g^i(s, x^i, v^i)}{\partial s} = g^i(s, x^i, v^i) \left[ \lambda_1 v^i + \frac{\lambda\lambda_2}{I} \sum_{i=1}^{I} \psi(x^i - x^j)(v^i - v^j) \right],$$

$$\frac{\partial g^i(s, x^i, v^i)}{\partial x^i} = g^i(s, x^i, v^i) \frac{\lambda\lambda_2}{I} \psi(v^i - v^j),$$

$$\frac{\partial g^i(s, x^i, v^i)}{\partial v^i} = g^i(s, x^i, v^i) \left[ s\lambda_1 + \frac{s\lambda\lambda_2}{I} \psi(x^i - x^j) \right],$$

$$\frac{\partial^2 g^i(s, x^i, v^i)}{\partial (x^i)^2} = g^i(s, x^i, v^i) \left[ \frac{s\lambda\lambda_2}{I} \psi(v^i - v^j) \right]^2,$$

$$\frac{\partial^2 g^i(s, x^i, v^i)}{\partial (v^i)^2} = g^i(s, x^i, v^i) \left[ s\lambda_1 + \frac{s\lambda\lambda_2}{I} \psi(x^i - x^j) \right],$$

and

$$\frac{\partial^2 g^i(s, x^i, v^i)}{\partial x^i \partial v^i} = g^i(s, x^i, v^i) \frac{\lambda\lambda_2}{I} \psi \left[ 1 + s\lambda_1 + \frac{s\lambda\lambda_2}{I} \psi(x^i - x^j) \right].$$

Therefore,

$$f^i(s, x^i, v^i, u^i) = \sum_{i=1}^{I} \exp(-\rho^i s) \alpha^i H_{01}^{i*} x^i v^i (u^i)^2 + g^i(s, x^i, v^i)$$
$$+ g^i(s, x^i, v^i) \left[ \lambda_1 v^i \right.$$
$$\left. + \frac{\lambda\lambda_2}{I} \sum_{i=1}^{I} \psi(x^i - x^j)(v^i - v^j) \right]$$
$$+ g^i(s, x^i, v^i) \frac{\lambda\lambda_2}{I} \psi u^i v^i (v^i - v^j)$$
$$+ g^i(s, x^i, v^i) \left[ s\lambda_1 + \frac{s\lambda\lambda_2}{I} \psi(x^i - x^j) \right]$$
$$\left[ \frac{\lambda}{I} \sum_{i=1}^{I} \psi |x^i - x^j| (v^i - v^j) \right]$$
$$+ \frac{1}{2} \left\{ (\sigma_1^i)^2 g^i(s, x^i, v^i) \left[ \frac{s\lambda\lambda_2}{I} \psi(v^i - v^j) \right] \right.$$
$$+ 2\rho(\sigma_1^i)^3 \frac{\lambda\lambda_2}{I} \psi g^i(s, x^i, v^i)$$
$$\times \left[ 1 + s\lambda_1 + \frac{s\lambda\lambda_2}{I} \psi(x^i - x^j) \right]$$
$$\left. + \sigma_2^{i*} g^i(s, x^i, v^i) \left[ s\lambda_1 + \frac{s\lambda\lambda_2}{I} \psi(x^i - x^j) \right]^2 \right\}. \tag{31}$$

In order to satisfy Eq. (14) either $\frac{\partial f^i(s, x^i, v^i, u^i)}{\partial u^i}$ or $\Psi_s^{i\tau}(x^i, v^i)$ has to be zero. As $\Psi_s^{i\tau}(x^i, v^i)$ is a wave function, it cannot be zero. Hence,

$$\frac{\partial f^i(s, x^i, v^i, u^i)}{\partial u^i} = 0,$$

or,

$$2u^i \exp(-\rho^i s) \alpha^i H_{01}^{i*} x^i v^i + g^i(s, x^i, v^i) \frac{\lambda\lambda_2}{I} \psi v^i (v^i - v^j)$$
$$+ g^i(s, x^i, v^i) \left[ s\lambda_1 + \frac{s\lambda\lambda_2}{I} \psi(x^i - x^j) \right]$$
$$\left[ \frac{\lambda}{I} \sum_{i=1}^{I} \psi |x^i - x^j| (v^i - v^j) \right] = 0.$$

Therefore, an optimal strategy of $i^{th}$ fish at time $s$ is

$$u^{i*}(s) = \frac{\exp\left\{\rho^i s + s\lambda_1 v^i + \frac{s\lambda\lambda_2}{I}\sum_{i=1}^{I}\psi(x^i - x^j)(v^i - v^i) + \lambda_3\sqrt{8/3}k^i(l)\right\}}{2\alpha^i H_{01}^{i*}x^i v^i}$$

$$\times \left\{\frac{\lambda\lambda_2}{I}v^i(v^j - v^i) + \left[s\lambda_1 + \frac{s\lambda\lambda_2}{I}\psi(x^i - x^j)\right]\right.$$

$$\left.\left[\frac{\lambda}{I}\sum_{i=1}^{I}\psi|x^i - x^j|(v^i - v^j)\right]\right\}, \tag{32}$$

such that $\alpha^i H_{01}^{i*}x^i v^i \neq 0$.

Suppose at time $s$ total number of fishes in a school is $I$ and $u^{i*}(s)$ be $i^{th}$ fish's strategy whether to stay in the school such that if its value is high, the fish will leave the school and vice versa. Assume communication rate $\psi$, nonnegative coupling strength $\lambda$, weight $\alpha^i$, Lagrangian multipliers $\{\lambda_1, \lambda_2, \lambda_3\}$ and discount rate $\rho^i$ take constant values. In the following cases I will discuss how other factors apart from the constants in Eq. (32) affect optimal strategy of fish $i$.

*Case I* Let the survival function $H_{01}^{i*} \in [0, 1]$ takes the value very close to zero because of some predatorial attacks which leads $u^{i*}(s)$ in Eq. (32) to take a very large value. Therefore, fish $i$ has to leave the school. Intuitively, after getting attacked by a predator fish $i$ is injured. Therefore, it cannot keep up with the speed of the school which lead it to decide to leave.

*Case II* Suppose fish $i$ is very close to fish $j$ in the school or $(x^i - x^j) \to 0$. Therefore,

$$u^{i*}(s) \to \frac{\exp\{\rho^i s + s\lambda_1 v^i + \lambda_3\sqrt{8/3}k^i(l)\}}{2\alpha^i H_{01}^{i*}x^i v^i}$$

$$\times \frac{\lambda\lambda_2}{I}v^i(v^j - v^i),$$

where $2\alpha^i I H_{01}^{i*}x^i v^i \neq 0$. As this value less than $u^{i*}(s)$ obtained in Eq. (32), fish $i$ will not leave the school. Intuitively, as any two fishes are very close to each other, they can withstand any external adverse environmental condition including attacks from predators. Hence, fish $i$'s strategy should be to stay in the school.

*Case III* If the relative velocity of fish $i$ and $j$ is similar or $(v^j - v^i) \to 0$, then

$$u^{i*}(s) \to \frac{\exp\{\rho^i s + s\lambda_1 v^i + \lambda_3\sqrt{8/3}k^i(l)\}}{2\alpha^i H_{01}^{i*}x^i v^i},$$

where $2\alpha^i H_{01}^{i*}x^i v^i \neq 0$. As the value of the optimal strategy gets reduced compared to Eq. (32), fish $i$ will not leave the school. Intuitively, as $v^i \to v^j$ for all $i \neq j$, all the fishes in the school have same relative velocity. Again similar to the argument as in **Case II**, they can survive any external attack, and fish $i$'s strategy is not to leave the school. Same result is obtained if the fish school is large or $I \to \infty$.

*Case IV* In Eq. (32), $k^i(l)$ represents a variant of a Gaussian free field (GFF and some harmonic function), where $l$ is some number coming from Brownian two-sphere. If $k^i(l)$ takes very high value, then optimal strategy $u^{i*}(s)$ goes up which leads to fish $i$ to leave the school. Intuitively, for a high value of $k^i$ would help increase the ergodicity of the strategy space due to external environmental conditions. This leads to break the fish school. Therefore, fish $i$'s optimal strategy is to leave.

## Conclusion

A Feynman type of path integral under $\sqrt{8/3}$-LQG surface has been introduced in this context of fish migration. The advantage of having Feynman type path integral is that it can be used in the case of generalized nonlinear stochastic differential equations where constructing Hamiltonian–Jacobi–Bellman (HJB) equation is impossible (Baaquie 2007). A few attempts of Feynman path integral have been used in order to determine animal behaviors (Kappen 2005b, 2007). Pontryagin maximum principle using HJB equation has been used rigorously in this literature. Furthermore, in the literature of fish migration, less number of contributions through stochastic differential equations have been made and all works are related to Pontryagin maximum principle. This paper is the first attempt where optimal strategy of a fish is determined through a new Feynman-type path integral approach where relative velocity of fish is on $\sqrt{8/3}$-LQG surface. When this surface takes the value of $\sqrt{8/3}$, it glues to a Brownian surface (Sheffield 2007; Duplantier and Sheffield 2011; Miller and Sheffield 2016b, a; Sheffield et al. 2016). The advantage of this method is that instead of using the properties of Brownian surface, one can replace it by a smooth function to do metric gluing.

In Proposition 1 a more generalized objective function than Yoshioka (2019) has been used subject to two stochastic dynamics expressed in Eqs. (2) and (3). Then a new Liouville–Feynman type of path integral method is constructed to get a Wick-rotated Schrödinger-type equation and the first-order condition with respect $u^i$ gives an optimal strategy of the $i^{th}$ fish to reach habitat $H_1$ through the migration process. In Example 1 an objective function similar to Yoshioka (2019), two stochastic differential equations similar to Uchitane et al. (2012); Nguyen et al. (2016) and Ton et al. (2014) have been used in order to determine an exact expression of $i^{th}$ fish's optimal strategy. In future research a more generalized tensor field would be used to find out an optimal strategy where actions of a fish falls under $p$-brane (Pramanik and Polansky 2019).

## Declarations

## References

Ahn SM, Ha SY (2010) Stochastic flocking dynamics of the cucker-smale model with multiplicative white noises. J Math Phys 51(10):103301

Arai T, Hayano H, Asami H, Miyazaki N (2003) Coexistence of anadromous and lacustrine life histories of the shirauo, salangichthys microdon. Fish Oceanogr 12(2):134–139

Arai T, Yang J, Miyazaki N (2006) Migration flexibility between freshwater and marine habitats of the pond smelt hypomesus nipponensis. J Fish Biol 68(5):1388–1398

Baaquie BE (1997) A path integral approach to option pricing with stochastic volatility: some exact results. J de Phys I 7(12):1733–1753

Baaquie BE (2007) Quantum finance: path integrals and Hamiltonians for options and interest rates. Cambridge University Press, Cambridge

Bauer S, Klaassen M (2013) Mechanistic models of animal migration behaviour-their diversity, structure and use. J Anim Ecol 82(3):498–508

Bochner S, Chandrasekharan K et al (1949) Fourier transforms. Princeton University Press, Princeton

Carrillo JA, Fornasier M, Rosado J, Toscani G (2010) Asymptotic flocking dynamics for the kinetic cucker-smale model. SIAM J Math Anal 42(1):218–236

Dorst JP (2019) Migration. https://www.britannica.com/science/migration-animal, [Online; posted 06-August-2019]

Duplantier B, Sheffield S (2011) Liouville quantum gravity and kpz. Inventiones mathematicae 185(2):333–393

Feynman RP (1948) Space-time approach to non-relativistic quantum mechanics. Rev Mod Phys 20(2):367

Fujiwara D (2017) Feynman's idea. In rigorous time slicing approach to feynman path integrals. Springer, Berlin

Guse B, Kail J, Radinger J, Schröder M, Kiesel J, Hering D, Wolter C, Fohrer N (2015) Eco-hydrologic model cascades: simulating land use and climate change impacts on hydrology, hydraulics and habitats for fish and macroinvertebrates. Sci Total Environ 533:542–556

Gwynne E, Miller J (2016) Metric gluing of brownian and $\sqrt{8/3}$-liouville quantum gravity surfaces. arXiv preprint arXiv:160800955

Ha SY, Lee K, Levy D et al (2009) Emergence of time-asymptotic flocking in a stochastic cucker-smale system. Commun Math Sci 7(2):453–469

Hebb DO (1949) The organization of behavior: a neuropsychological theory. Wiley, Chapman & Hall

Hellman Z, Levy YJ (2019) Measurable selection for purely atomic games. Econometrica 87(2):593–629

Jonsson N, Jonsson B (2002) Migration of anadromous brown trout salmo trutta in a Norwegian river. Freshw Biol 47(8):1391–1401

Kappen HJ (2005a) Linear theory for control of nonlinear stochastic systems. Phys Rev Lett 95(20):200201

Kappen HJ (2005b) Path integrals and symmetry breaking for optimal control theory. J Stat Mech Theor Exp 2005(11):P11011

Kappen HJ (2007) An introduction to stochastic control theory, path integrals and reinforcement learning. AIP conference proceedings, American Institute of Physics 887:149–181

Knizhnik VG, Polyakov AM, Zamolodchikov AB (1988) Fractal structure of 2d-quantum gravity. Modern Phys Lett A 3(08):819–826

Lande R, Engen S, Saether BE et al (2003) Stochastic population dynamics in ecology and conservation. Oxford University Press, Oxford

Larsson M (2012) Incidental sounds of locomotion in animal cognition. Anim Cognit 15(1):1–13

Lemasson BH, Haefner JW, Bowen MD (2014) Schooling increases risk exposure for fish navigating past artificial barriers. PloS One 9(9):108220

Marcet A, Marimon R (2019) Recursive contracts. Econometrica 87(5):1589–1631

Mas-Colell A, Whinston MD, Green JR et al (1995) Microeconomic theory, vol 1. Oxford University Press, New York

Miller J (2018) Liouville quantum gravity as a metric space and a scaling limit. In: proceedings of the international congress of mathematicians: Rio de Janeiro 2018, World Scientific, pp 2945–2971

Miller J, Sheffield S (2016a) Imaginary geometry i: interacting sles. Prob Theory Relat Fields 164(3–4):553–705

Miller J, Sheffield S (2016b) Liouville quantum gravity and the brownian map iii: the conformal structure is determined. arXiv preprint arXiv:160805391

Nguyen LTH, Ta VT, Yagi A (2016) Obstacle avoiding patterns and cohesiveness of fish school. J Theor Biol 406:116–123

Øksendal B (2003) Stochastic differential equations. In: stochastic differential equations. Springer, Berlin

Øksendal B, Sulem A (2019) Applied Stochastic Control of Jump Diffusions. Springer. This a later edition of the original book published in 2007. Øksendal, B, & Sulem A. (2007). Applied Stochastic Control of Jump Diffusions (Vol. 498). Berlin: Springer. https://link.springer.com/content/pdf/10.1007/978-3-030-02781-0.pdf

Pinti J, Celani A, Thygesen UH, Mariani P (2020) Optimal navigation and behavioural traits in oceanic migrations. Theor Ecol 13(4):583–593

Pitici M (2018) The best writing on mathematics 2017. Princeton University Press, Princeton

Polyakov AM (1981) Quantum geometry of bosonic strings. Phys Lett B 103(3):207–210

Polyakov AM (1987) Quantum gravity in two dimensions. Modern Phys Lett A 2(11):893–898

Polyakov AM (1996) Quantum geometry of fermionic strings. In: 30 Years Of The Landau Institute—Selected Papers, World Scientific, pp 602–604

Pramanik P (2020) Optimization of market stochastic dynamics. SN Op Res Forum, Springer 1:1–17

Pramanik P, Polansky AM (2019) Semicooperation under curved strategy spacetime. arXiv preprint arXiv:191212146

Pramanik P, Polansky AM (2020a) Motivation to run in one-day cricket. arXiv preprint arXiv:200111099

Pramanik P, Polansky AM (2020b) Optimization of a dynamic profit function using euclidean path integral. arXiv preprint arXiv:200209394

Radinger J, Wolter C (2015) Disentangling the effects of habitat suitability, dispersal, and fragmentation on the distribution of river fishes. Ecol Appl 25(4):914–927

Ross K (2008) Stochastic control in continuous time. Lecture Notes on Continuous Time Stochastic Control, Spring

Schramm O (2000) Scaling limits of loop-erased random walks and uniform spanning trees. Israel J Math 118(1):221–288

Shang Y, Bouffanais R (2014) Influence of the number of topologically interacting neighbors on swarm dynamics. Sci Rep 4:4184

Sheffield S (2007) Gaussian free fields for mathematicians. Prob Theory Relat Fields 139(3–4):521–541

Sheffield S et al (2016) Conformal weldings of random surfaces: Sle and the quantum gravity zipper. Ann Prob 44(5):3474–3545

Simon B (1979) Functional integration and quantum physics. Academic press, Cambridge

Theodorou E, Buchli J, Schaal S (2010) Reinforcement learning of motor skills in high dimensions: a path integral approach. In: Robotics and Automation (ICRA), 2010 IEEE International Conference on, IEEE, pp 2397–2403

Theodorou EA (2011) Iterative path integral stochastic optimal control: theory and applications to motor control. University of Southern California. http://citeseerx.ist.psu.edu/viewdoc/download?doi=10.1.1.418.8228&rep=rep1&type=pdf

Ton TV, Linh NTH, Yagi A (2014) Flocking and non-flocking behavior in a stochastic cucker-smale system. Anal Appl 12(01):63–73

Torney CJ, Lorenzi T, Couzin ID, Levin SA (2015) Social information use and the evolution of unresponsiveness in collective systems. J Royal Soc Interf 12(103):20140893

Uchitane T, Ton TV, Yagi A (2012) An ordinary differential equation model for fish schooling. Sci Math Jpn 75(3):339–350

Van Den Broek B, Wiegerinck W, Kappen B (2008) Graphical model inference in optimal control of stochastic multi-agent systems. J Art Intell Res 32:95–122

Wang X, Pan Q, Kang Y, He M (2016) Predator group size distributions in predator-prey systems. Ecol complex 26:117–127

Yang I, Morzfeld M, Tomlin CJ, Chorin AJ (2014) Path integral formulation of stochastic optimal control with generalized costs. IFAC proceedings volumes 47(3):6994–7000

Yoshioka H (2017) A simple game-theoretic model for upstream fish migration. Theory Biosci 136(3–4):99–111

Yoshioka H (2019) A stochastic differential game approach toward animal migration. Theory Biosci 138(2):277–303

Yoshioka H, Yaegashi Y (2018) An optimal stopping approach for onset of fish migration. Theory Biosci 137(2):99–116

Yoshioka H, Shirai T, Tagami D (2019) A mixed optimal control approach for upstream fish migration. J Sustain Develop Energy W Environ Syst 7(1):101–121